\documentclass{amsart}

\usepackage{amssymb}
\usepackage{graphicx}
\usepackage{latexsym}
\usepackage{epsfig}
\usepackage{psfrag}
\def\scfig #1 #2 {\resizebox{#2}{!}{\includegraphics{#1}}}
\def\C{\mathbb C}
\def\CE{\mathcal E}
\def\qed{\hfill$\Box$}
\def\bi{\bibitem}

\newtheorem{theorem}{Theorem}
\newtheorem{lemma}[theorem]{Lemma}
\newtheorem{proposition}[theorem]{Proposition}
\newtheorem{corollary}[theorem]{Corollary}

\theoremstyle{definition}
\newtheorem{definition}[theorem]{Definition}

\theoremstyle{remark}
\newtheorem{remark}[theorem]{Remark}

\newtheorem{Theorem}{\sc Theorem}

\def\hpic #1 #2 {\mbox{$\begin{array}[c]{l} 
\epsfig{file=#1,height=#2}\end{array}$}}
\def\wpic #1 #2 {\mbox{$\begin{array}[c]{l} 
\epsfig{file=#1,width=#2}\end{array}$}}

\def\C{\mathbb C}

\def\Z{\mathbb Z}

\def\CE{{\mathcal {E}}}

\def\CH{{\mathcal {H}}}

\def\CL{{\mathcal {L}}}

\def\be{\begin{equation}}
\def\ee{\end{equation}}
\def\ba{\begin{array}}
\def\ea{\end{array}}
\def\bt{\begin{Theorem}}
\def\et{\end{Theorem}}
\def\bi{\begin{itemize}}
\def\ei{\end{itemize}}
\def\bea{\begin{eqnarray}}
\def\eea{\end{eqnarray}}
\def\beast{\begin{eqnarray*}}
\def\eeast{\end{eqnarray*}}
\def\ben{\begin{enumerate}}
\def\een{\end{enumerate}}

\def\qed{\hfill$\Box$}
\def\bi{\bibitem}

\def\rar{\rightarrow}
\def\Rar{\Rightarrow}

\def\a{{\overline{a}}}

\def\u{{\overline{u}}}

\begin{document}
\baselineskip 16pt

\title[The GJS-construction for graphs]{On the Guionnet-Jones-Shlyakhtenko construction for graphs}

\author{Vijay Kodiyalam}
\address{The Institute of Mathematical Sciences, Chennai, India}
\email{vijay@imsc.res.in}

\author{V. S. Sunder}
\address{The Institute of Mathematical Sciences, Chennai, India}
\email{sunder@imsc.res.in}





\begin{abstract}
Using an analogue of the Guionnet-Jones-Shlaykhtenko construction
for graphs we show that their construction
applied to  any subfactor planar algebra of finite depth yields
 an inclusion of  interpolated
free group factors with finite parameter, thereby giving another proof of their universality
 for finite depth planar algebras.
\end{abstract}

\maketitle



The main theorem of \cite{GnnJnsShl2008} constructs an extremal finite index $II_1$ subfactor $N = M_0 \subseteq M_1 = M$ from a subfactor planar
algebra $P$ with the property that the planar algebra of $N \subseteq M$ is isomorphic to $P$. We show in this paper
that if $P$ is a
subfactor planar algebra of modulus $\delta>1$ and of finite depth, then, for the associated subfactor $N \subseteq M$, 
there are isomorphisms
$N \cong LF(r)$ and $M \cong LF(s)$  for some $1 < r,s < \infty$,
where $LF(t)$ for $1 <t < \infty$ is the interpolated free group factor 
of \cite{Dyk1994} and \cite{Rdl1994}.
This can be regarded as yet another proof of the fact -
see \cite{Rdl1994} and \cite{Dyk2002} -  that interpolated free group factors with finite parameter are universal for finite depth subfactor
planar algebras. The word `universal' above is used in the sense of \cite{PpaShl2003} where they essentially prove 
that $LF(\infty)$ is universal for all subfactor planar algebras.

We shall now outline the structure of this paper. In \S 1 we construct - see Proposition \ref{p1} - 
a graded, tracial, faithful $*$-probability space $Gr(\Gamma)$  associated
to a  finite, weighted, bipartite graph $\Gamma$ and establish - see Proposition \ref{same} - an isomorphism between $Gr(\Gamma)$ and a filtered, tracial, faithful $*$-probability space $F(\Gamma)$ - see Proposition \ref{p3}. Our main interest will be in an associated finite von Neumann algebra $M(\Gamma)$ and  some of its corners determined by sets
of vertices of $\Gamma$ - specifically the corner $M(\Gamma,0)$ (respectively $M(\Gamma,1)$) determined by the set of even (respectively odd) vertices of $\Gamma$.
The main result in \S 2 
asserts - see Theorem \ref{Mgamstr} - 
that if $\Gamma$ is a connected graph with more than one edge,
then, $M(\Gamma)$ is the direct sum of a $II_1$ factor and a  finite-dimensional abelian algebra.
The goal of \S 3 is to 
express $Gr(\Gamma,0)$ and $M(\Gamma,0)$ - see Proposition \ref{freeness} and equation (\ref{mgamam}) - as amalgamated free
products of the corresponding algebras associated to subgraphs with
a single odd vertex.
In \S 4 we determine the structure of $M(\Lambda,0)$ - see Corollary \ref{l0str} - for a graph
$\Lambda$ with a single odd vertex.
The penultimate \S 5 proves - see Theorem \ref{mgamlf} - one of our main results : for a connected graph $\Gamma$ with more than one edge and equipped with its
Perron-Frobenius weighting, the algebra $M(\Gamma)$ is an 
interpolated free group factor with finite parameter.
The final \S 6 applies this - see Theorems \ref{main} and \ref{main2} -  to show that the Guionnet-Jones-Shlyakhtenko (henceforth GJS) construction applied
to a finite depth subfactor planar algebra yields an inclusion of interpolated free group factors with finite parameters.

\section{The global graded probability space associated to a graph}

The goal of this section is to associate
a
graded, tracial, faithful $*$-probability space $Gr(\Gamma)$  
and a von Neumann algebra $M(\Gamma)$
to a graph $\Gamma$ . Recall that a tracial $*$-probability space
consists of a unital, complex $*$-algebra $A$ equipped with a trace
$\tau : A \rightarrow {\mathbb C}$ that satisfies $\tau(1) = 1$
 and  $\tau(a^*a) \geq 0$, for
all $a \in A$. It is said to be graded if the algebra $A$ is graded and
to be faithful if $\tau(a^*a) = 0 \Rightarrow a=0$.

Throughout this paper, by a graph, we will mean a finite, weighted, bipartite graph which consists of the following data: 
(i) a finite set $V$ of `vertices' partitioned as $V_0
\coprod V_1$ - the sets $V_0$ and $V_1$ will be referred to as
sets of even and odd vertices respectively, (ii) a finite set $E$ of `edges'
equipped with `start' and `finish' maps $s,f : E \rightarrow V$ and
a `reversal' involution $\xi \mapsto \tilde{\xi}$ of $E$ intertwining $s$ and
$f$ such that $s(\xi) \in V_0 \Leftrightarrow f(\xi) \in V_1$, and 
(iii) a `weighting' which is a  function $\mu : V \rightarrow {\mathbb R}_+$ normalised such that
$\sum_{v \in V} \mu^2(v) = 1$.

For us, the main examples of such graphs are the principal graphs of non-trivial $II_1$-subfactors of finite depth (where 
$\mu$ is given by the square root of
an appropriately normalised Perron-Frobenius eigenvector)
and their subgraphs (with 
the restricted $\mu$ appropriately normalised).

The construction of $Gr(\Gamma)$
 involves paths in $\Gamma$, notations and definitions for which we discuss briefly.
A path $\xi$ in $\Gamma$ is denoted
$$
(v^\xi_0 \stackrel{\xi_1}{\rightarrow} v^\xi_1 \stackrel{\xi_2}{\rightarrow}
v^\xi_2 \stackrel{\xi_3}{\rightarrow} \cdots \stackrel{\xi_n}{\rightarrow}
v^\xi_n),
$$
where $v^\xi_i \in V$ and $(v^\xi_{i-1} \stackrel{\xi_i}{\rightarrow} v^\xi_i) \in E$, with the notation being self-explanatory.
The start and finish vertex functions on paths in $\Gamma$ will also be
denoted by $s(\cdot)$ and $f(\cdot)$ respectively and the length function by $\ell(\cdot)$, so that $s(\xi) = v^\xi_0, f(\xi) = v^\xi_n$
and $\ell(\xi) = n$.
For  $0 \leq i \leq j \leq n$, we will use
notation such as $\xi_{[i,j]}$ for the path $(v^\xi_i \stackrel{\xi_{i+1}}{\rightarrow}  v^\xi_{i+1} \stackrel{\xi_{i+2}}{\rightarrow} v^\xi_{i+2} \cdots 
\stackrel{\xi_j}{\rightarrow} v^\xi_j)$, where the interval refers to the
vertex indices.
The symbol $\circ$ will denote composition of paths and $\widetilde{\ \ }$ will stand for path reversal.
For $n \geq 0$, the path space $P_n(\Gamma)$ associated to the graph $\Gamma$
is the complex vector space with basis $\{ [\xi]: \xi$ is a path of length
$n$ in $\Gamma\}$.

We will now define  $Gr(\Gamma)$ and its structure maps. As a graded
vector space, $Gr(\Gamma) = \oplus_{n \geq 0} P_n(\Gamma)$.
The multiplication in $Gr(\Gamma)$, denoted by $\bullet$, is given
by concatenation on the path basis and extended by linearity:
$$
[\xi] \bullet [\eta] = 
\left\{ 
\begin{array}{ll}
                   0 & {\text {if\ }} f(\xi) \neq s(\eta) \\
                   {[}\xi \circ \eta{]} & {\text {if\ }} f(\xi) = s(\eta).
\end{array}
\right. 
$$
The involution $*$ on $Gr(\Gamma)$ is defined by conjugate
linear extension of the reversal map $\widetilde{\ \ }$ on the path basis, i.e.,
$[\xi]^* = [\tilde{\xi}]$. 
We define a linear functional $\tau$
on $Gr(\Gamma)$ motivated by the GJS trace.
Suppose that $[\xi] \in P_n(\Gamma)$ for $n \geq 1$. Define 
$$
\tau([\xi]) =  \sum_T \tau_T([\xi])
$$
where the sum is over all Temperley-Lieb equivalence
relations\footnote{These are the non-crossing relations with every
  class having two elements.}
$T$ on $\{1,2,\cdots,n\}$ (so that it is an empty sum, hence vanishes, for
$n$ odd) and $\tau_T$ is defined
by
$$
\tau_T([\xi]) = 
\prod_{\{\{i,j\} \in T : i <j\}} \delta_{\xi_i,\widetilde{\xi_j}}
\prod_{C \in K(T)} \mu(v^\xi_C)^{2-|C|}
$$
where (i) $K(T)$ is the Kreweras complement of $T$ - see \cite{NcaSpc2006} - which is also 
a non-crossing partition of $\{1,2,\cdots,n\}$ and (ii) $v^\xi_C =
v^\xi_c$ for any $c \in C$ (all of which must be equal if the
first product is non-zero). When $n=0$,  we set $\tau([(v)]) = \mu^2(v)$. 
 
\begin{proposition}\label{p1}
$Gr(\Gamma)$ is a graded, unital, associative, $*$-algebra and $\tau$
is a normalised trace on $Gr(\Gamma)$.
\end{proposition}

\begin{proof}
The only not completely obvious assertion is the traciality
of $\tau$, which too follows, after a little thought, from the rotational invariance of the set
of all TL-equivalence relations and from the definition of the product in $Gr(\Gamma)$.
\end{proof}

Note that the multiplicative identity of $Gr(\Gamma)$ is the element $\sum_{v \in V} [(v)] \in P_0(\Gamma)$.
In view of the fact that the different $[(v)]$, for $v \in V$, are orthogonal idempotents (adding to $1$), we will denote $[(v)]$ also by
$e_v$. 
It is useful to observe that an element, say $x$, of $P_n(\Gamma)$
may be regarded as the square matrix, with rows and
columns indexed by $V$, with $(v,w)$ entry given by
$e_vxe_w$ (the part of $x$ which is a linear combination
of paths beginning at $v$ and ending at $w$).

The proof of positivity and faithfulness of the trace $\tau$ involves
some work with a different avatar of $Gr(\Gamma)$ which we will
find very useful. We begin by recalling, from \cite{JnsShlWlk2008},
the category epi-TL which we will denote by ${\mathcal E}$. The objects of ${\mathcal E}$ are denoted $[n]$ for $n \geq 0$ and thought
of as $n$-points (labelled $1,2,3,\cdots,n$) arranged on a horizontal line. A morphism in $Hom([n],[m])$
consists of a rectangle with $m$-points on the top horizontal line,
$n$-points on the bottom horizontal line and a Temperley-Lieb like
tangle in between, subject to the restriction that each of the points
above is joined to a point below. It must be observed that
$Hom([n],[m])$ is non-empty precisely when $n-m$ is a non-negative even integer.
Morphisms are composed by vertical stacking.

The morphisms in ${\mathcal E}$ are generated by those which have a
single cap on the bottom line. 
Let $S^n_i : [n] \rightarrow [n-2]$ (for $1 \leq i < n$) denote the generator with the $i^{th}$ and $(i+1)^{st}$
points on the bottom line capped. Some work shows that all
relations among the morphisms are consequences of the relations
\be\label{snrel}
S^{n-2}_p S^n_q = S^{n-2}_{q}S^n_{p+2}
\ee
for  $n-2 > p \geq q \geq 1$.
In fact any element of $Hom([m+2k],[m])$
is uniquely expressible in the form $S^{m+2}_{i_1}S^{m+4}_{i_2}\cdots
S^{m+2k}_{i_k}$ with $1 \leq i_1 < i_2 < \cdots < i_k < m+2k$. (The
left end points of the $k$ caps of the morphism are precisely at the
places $i_1, i_2, \cdots , i_k$.) Such a morphism will be called 
  non-nested if the caps are `not nested', or equivalently, if
$i_{j+1} \geq i_j + 2$ for each $j < k$ in its `canonical
decomposition' as above.

It follows that the category ${\mathcal E}$ `acts' on the collection of vector
spaces $P_n(\Gamma)$ in the sense that any element of $Hom([n],[m])$ yields a vector space homomorphism $P_n(\Gamma) \rightarrow
P_m(\Gamma)$ with this assignment being compatible with
compositions on both sides. Such an action can be defined\footnote{Thus we are saying that the operators defined by equation (\ref{snact}) satisfy the relations (\ref{snrel}).} with
$S^n_i$ acting by
\be\label{snact}
S^n_i([\xi]) = \delta_{\xi_i,\widetilde{\xi_{i+1}}}
\frac{\mu(v^{\xi}_i)}{\mu(v^\xi_{i\pm 1})} [\xi_{[0,i-1]} \circ
\xi_{[i+1,n]}]
\ee
for $[\xi] \in P_n(\Gamma)$.
More generally, given an arbitrary $S \in Hom([n],[m])$, it specifies
a partition of $[n]$ as $T \cup E$, where $T$ is the subset of
points in $[n]$ that are joined to a point in $[m]$ and $E$ is its
complement. It also specifies a Temperley-Lieb  equivalence relation $\sim$ on $E$. The action of $S$ is then explicitly given by
\begin{equation}\label{sform}
S([\xi]) =  \prod_{\{i,j\} \in\  \sim: i < j} \left( \delta_{\xi_i,\widetilde{\xi_j}}
\frac{\mu(v^{\xi}_i)}{\mu(v^\xi_{j})} \right) [\circ_{t \in T} \xi_t],
\end{equation}
where the concatenation is done in increasing order of elements of $T$
and is interpreted as $[(f(\xi))]$ if $T = \emptyset.$ (As in the
equations displayed above, we shall often identify elements of
$Hom([n],[m])$ with the associated operators from $P_n(\Gamma)$ to
$P_m(\Gamma)$.)

The following lemma is a special case (of Proposition \ref{diffexp}) which  both motivates and is used in the proof
of  a different expression  for
$S([\xi])$ when $S \in Hom([2n],[0])$. 
Note that in this case, $E=\{1,2,\cdots,2n\}$ and $\sim \ = S$ regarded as an equivalence relation.

\begin{lemma}\label{diffexplem}
Let $[\xi] \in P_{2n}(\Gamma)$ and $S \in Hom([2n],[0])$ be
given by $\{\{1,2n\},\{2,2n-1\},\cdots,\{n,n+1\}\}$. Then,
\begin{eqnarray*}
S([\xi]) = \frac{\mu(v^\xi_n)}{\mu(v^\xi_{2n})}
\prod_{i=1}^n \delta_{\xi_i,\widetilde{\xi_{2n+1-i}}}
  \times [(v^\xi_{2n})]. \nonumber \\
\end{eqnarray*}
\end{lemma}

\begin{proof}
We may assume that $\xi$ is a path consistent with $S$
in the sense that $\xi_i = \widetilde{\xi_j}$ whenever $\{i,j\} \in S$,
since otherwise, both sides of the desired equality vanish.
Thus, $\xi_i = \widetilde{\xi_{2n+1-i}}$ and in particular,  $v^\xi_i = v^\xi_{2n-i}$ for each $i = 0,1,\cdots,2n$.

Using equation (\ref{sform}), it now suffices to check that
$$
\prod_{i=1}^n \frac{\mu(v^\xi_i)}{\mu(v^\xi_{2n+1-i})} = \frac{\mu(v^\xi_n)}{\mu(v^\xi_{2n})}.
$$
But substituting $v^\xi_i = v^\xi_{2n-i}$, we see that the product on the left telescopes to the expression on the right.
\end{proof}

We next treat the case of a general $S$.

\begin{proposition}\label{diffexp}
For $[\xi] \in P_{2n}(\Gamma)$ and $S \in Hom([2n],[0])$,
\begin{eqnarray}
S([\xi]) = \frac{\mu(v^\xi_n)}{\mu(v^\xi_{2n})}
  \left( \prod_{\{i,j\} \in S: i < j \leq n} \delta_{\xi_i,\widetilde{\xi_j}}
\frac{\mu(v^{\xi}_i)}{\mu(v^\xi_{j})} \right)
\left( \prod_{\{i,j\} \in S: i \leq n < j} \delta_{\xi_i,\widetilde{\xi_j}}
 \right) \times \nonumber \\
   \left( \prod_{\{i,j\} \in S: n < i < j} \delta_{\xi_i,\widetilde{\xi_j}}
\frac{\mu(v^{\xi}_i)}{\mu(v^\xi_{j})} \right) [(v^\xi_{2n})]. \label{sndef}
\end{eqnarray}
\end{proposition}

\begin{proof}[Sketch of Proof] As in the proof of Lemma \ref{diffexplem},
we may assume that $\xi$ is a path consistent with $S$.
In this case,
comparison with equation (\ref{sform})
now shows that it suffices to see the following:
\begin{equation}\label{desired}
\prod_{\{i,j\} \in S: i < n < j} \left(
\frac{\mu(v^{\xi}_i)}{\mu(v^\xi_{j})} \right) = \frac{\mu(v^\xi_n)}{\mu(v^\xi_{2n})}.
\end{equation}

We illustrate by way of an example why this holds. Consider the
$S$ in Figure \ref{eqrel2} which corresponds to the equivalence
relation $\{\{1,10\},\{2,7\},\{3,6\},\{4,5\},\{8,9\}\}$
\begin{figure}[htbp]
\begin{center}
\includegraphics[height=2cm]{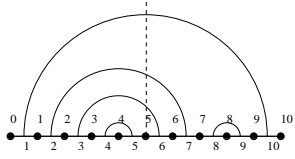}
\caption{The element $S \in Hom([10],[0])$}
\label{eqrel2}
\end{center}
\end{figure}
The numbers $1,2,\cdots,10$ below the line index the edges of $\xi$
while the numbers $0,1,\cdots,10$ above index the vertices of $\xi$.
The LHS of equation (\ref{desired}) in this example is
$$
\frac{\mu(v^\xi_1)}{\mu(v^\xi_{10})} \frac{\mu(v^\xi_2)}{\mu(v^\xi_{7})} \frac{\mu(v^\xi_3)}{\mu(v^\xi_{6})}.
$$
The point now is that when the Kronecker delta terms are all non-zero,
all the $v^\xi_i$ in a single `region' are equal. Thus in this example,
$v^\xi_0 = v^\xi_{10}$, $v^\xi_1 =v^\xi_7 = v^\xi_9$, $v^\xi_2 = v^\xi_6$ and
$v^\xi_3 = v^\xi_5$. Hence, after cancellation, the LHS does simplify
to the RHS. 

Even in general, it should be clear that this happens. For the LHS of
equation (\ref{desired}) does not depend on those classes $\{i,j\}$ of $S$ for which both $i,j$ are either (i)~at most $n$ or (ii)~at least $n+1$.
Observing that the numbers of classes satisfying (i) and (ii) are equal,  we delete these classes and then we are in a situation
where Lemma \ref{diffexplem} applies.
\end{proof}

We will next define the algebra $F(\Gamma)$ and its structure maps.
As a vector space, $F(\Gamma) = \oplus_{n \geq 0} P_n(\Gamma)$. The multiplication, denoted $\#$, is defined as follows
on the path basis and extended linearly. Given $[\xi]
\in P_m(\Gamma)$ and $[\eta] \in P_n(\Gamma)$, the
product $[\xi] \# [\eta]$ has a component in $P_{m+n-2k}(\Gamma)$
for $0 \leq k \leq min\{m,n\}$, this component being given by
\begin{eqnarray*}
\lefteqn{([\xi] \# [\eta])_{m+n-2k} =} \\
&\left\{
\ba{ll}
[\xi] \bullet [\eta] &\mbox {if } k=0\\
S_{m-k+1}^{m+n-2(k-1)}S^{m+n-2(k-2)}_{m-k+2}\cdots S^{m+n-2}_{m-1}S^{m+n}_m([\xi] \bullet [\eta]) &\mbox {if } k>0
\ea
\right.
\end{eqnarray*}
The $*$ on $F(\Gamma)$ is exactly the same as that on $Gr(\Gamma)$ - namely $[\xi]^* = [\tilde{\xi}]$ extended conjugate linearly.
Finally, define a linear functional $t$ on $F(\Gamma)$ by setting its
restriction to $P_n(\Gamma)$ for $n \geq 1$ to be $0$ and by
linearly extending the map $[(v)] \mapsto \mu^2(v)$ on $P_0(\Gamma)$.

\begin{proposition}\label{p3}
$F(\Gamma)$ is a unital, associative, $*$-algebra and $t$
is a faithful, positive trace on $F(\Gamma)$.
\end{proposition}

\begin{proof}
A proof very similar to that in \cite{KdySnd2008}, and which we consequently omit, shows that $F(\Gamma)$
is a unital, associative $*$-algebra.
To show that $t$ is a faithful, positive trace it suffices to check that $\langle x , y \rangle = t(y^*x)$
defines an inner-product on $F(\Gamma)$
satisfying $\langle x , y \rangle = \langle y^* , x^* \rangle$. Consider the path basis
$[\xi]$ of $F(\Gamma)$. It follows from the definitions and Lemma \ref{diffexplem} that
$\langle [\xi] , [\eta] \rangle = \delta_{\xi,\eta} \mu(s(\xi))\mu(f(\xi))$, finishing the proof.
\end{proof}

We next define maps $\phi : Gr(\Gamma) \rightarrow F(\Gamma)$
and $\psi : F(\Gamma) \rightarrow Gr(\Gamma)$ as follows. Each of
these restricts to maps from $P_n(\Gamma)$ to $\oplus_{m=0}^n
P_m(\Gamma)$. Consequently, the maps $\phi, \psi$ may be represented
by upper-triangular matrices $((\phi^m_n))$ and $((\psi^m_n))$ where
$\phi^m_n, ~\psi^m_n : P_n(\Gamma) \rightarrow P_m(\Gamma)$ are zero if
$m> n$. We define $\phi^m_n$ to be the (action by the) sum of all
elements of $Hom([n],[m])$ and $\psi^m_n$ to be $(-1)^{n-m}$ times the
(action by the) sum of all the non-nested elements of $Hom([n],[m])$.

We now have the following proposition that identifies 
$Gr(\Gamma)$ and $F(\Gamma)$.

\begin{proposition}\label{same}
The maps $\phi$ and $\psi$ define mutually inverse $*$-isomorphisms between $Gr(\Gamma)$ and $F(\Gamma)$
that intertwine the functionals $\tau$ and $t$.
\end{proposition}

The proof uses the following lemma about the Kreweras complement
of Temperley-Lieb equivalence relations.

\begin{lemma} 
Let $S$ be a Temperley-Lieb equivalence relation on $\{1,2,\cdots,2n\}$ and $K(S)$ be its Kreweras complement. Then, for any class
$C = \{a_1,\cdots,a_k\}$ of $K(S)$ with $a_1 < \cdots < a_k$,
all the $a_i$ have the same parity and $\{a_i+1,a_{i+1}\} \in S$
for each $i = 1,\cdots,k$ (where $a_i+1$ is computed modulo $2n$ and  $i+1$ is computed modulo $k$).
\end{lemma}

\begin{proof}
Induce on $n$, with the basis case $n=1$ following by a direct check.
For $n>1$ take $i \leq 2n-1$ largest so that $\{i,i+1\} \in S$.  Let $T = S|_{\{1,2,\cdots,2n\} \setminus \{i,i+1\}}$. The Kreweras complement of $S$ is obtained from that of $T$ by adding 
$i+1$ to the class of $i-1$ and adding the singleton class $\{i\}$. Observe that $i+1$ is the largest element in its $K(S)$ class by choice of $i$.
Now by induction, the parity assertion holds and further, the new $\{a_i+1,a_{i+1}\}$ that are needed to be shown to belong to $S$ are both $\{i,i+1\}$ which is, indeed, in $S$.
\end{proof}

\begin{proof}[Proof of Proposition \ref{same}]
The proof that the maps $\phi$ and $\psi$ define mutually inverse $*$-isomorphisms between $Gr(\Gamma)$ and $F(\Gamma)$ is nearly identical to that of Lemma 5.1 in \cite{JnsShlWlk2008}  and depends essentially
only on properties of the category ${\mathcal E}$. We omit it here.

The intertwining assertion that needs to be checked is that
$\tau = t \circ \phi$ on $Gr(\Gamma)$.
Note that both sides vanish on paths of odd length and that
if $[\xi]$ is a path of length $2n$, then,  $\tau([\xi]) = \sum_T \tau_T([\xi])$ where the sum is over all Temperley-Lieb
relations $T$ on $\{1,2,\cdots,2n\}$ while $t \circ \phi([\xi]) =t (\sum_S 
S([\xi]))$ where the sum is over all $S \in Hom([2n],[0])$,
since $t$ vanishes on paths of positive length.
The natural identification between Temperley-Lieb equivalence relations
on $\{1,2,\cdots,2n\}$ and $Hom([2n],[0])$ shows that it suffices to see
that $\tau_S([\xi]) = t \circ S([\xi])$ for any Temperley-Lieb relation $S$ on $\{1,2,\cdots,2n\}$.
Now both these vanish unless $s(\xi) = f(\xi)$;  so we assume this.
Unravelling the definitions, we need to see that under these assumptions,
$$
\prod_{\{\{i,j\} \in S : i <j\}} \delta_{\xi_i,\widetilde{\xi_j}}
\prod_{C \in K(S)} \mu(v^\xi_C)^{2-|C|}
=\mu(v^\xi_{2n})^{2} \prod_{\{\{i,j\} \in S : i < j\}} \left( \delta_{\xi_i,\widetilde{\xi_j}}
\frac{\mu(v^{\xi}_i)}{\mu(v^\xi_{j})} \right),
$$
with $K(S)$ being the Kreweras complement of $S$ and 
$v^\xi_C =v^\xi_c$ for any $c \in C$.
We may further assume that $\xi$ is a path consistent with $S$
in the sense that $\xi_i = \widetilde{\xi_j}$ whenever $\{i,j\} \in S$
and show the following
$$
\mu(v^\xi_{2n})^{-2} \prod_{C \in K(S)} \mu(v^\xi_C)^{2-|C|}
=\prod_{\{\{i,j\} \in S : i < j\}} \left( \frac{\mu(v^{\xi}_i)}{\mu(v^\xi_{j})} \right).
$$

The product on the right may be rewritten as $\prod_{i=1}^{2n} \mu(v^\xi_i)^{\epsilon_S(i)}$ where $\epsilon_S(i)$ is $1$ or $-1$ according
as $i$ is the smaller or larger element in its $S$-class.
Next, we may regroup this product in terms of classes of $K(S)$
as $\prod_{C \in K(S)} \prod_{c \in C} \mu(v^\xi_c)^{\epsilon_S(c)}$.
Now, as we have observed before, if $\xi$ is consistent with $S$ then
all $v^\xi_c$ for $c \in C$ are equal (to a vertex denoted $v^\xi_C$) and so this product now becomes
$\prod_{C \in K(S)} \mu(v^\xi_C)^{\sum_{c \in C} \epsilon_S(c)}$.
Comparing with the product on the left, what needs to be seen is that
if $C$ is any class of $K(S)$ then
\begin{equation}\label{epsilon}
\sum_{c \in C} \epsilon_S(c) =
\left\{
\begin{array}{ll}
2-|C| &{\mbox {if\ \ }} 2n \notin  C\\
-|C|    &{\mbox {if\ \ }} 2n \in C
\end{array}
\right.
\end{equation}

To prove equation (\ref{epsilon}), it suffices to see 
that for any non-external (i.e., not containing $2n$) class $C$
of $K(S)$, $\epsilon_S(c)$ is $1$ or $-1$ according as $c$ is the smallest
element in $C$ or not, while for the external class, $\epsilon_S(c) = -1$
for all its elements.
But this is an easy consequence of Lemma \ref{epsilon}.
If $C = \{a_1,\cdots,a_k\}$ is a $K(S)$ class for which $a_k \neq 2n$, the definition
of $\epsilon_S$ (together with Lemma \ref{epsilon}) shows that $\epsilon_S(a_1) = 1$ while $\epsilon_S(a_i) 
= -1$ for $i \geq 2$. On the other hand, if $a_k = 2n$, then it similarly follows that all $\epsilon_S(a_i) =-1$, completing the proof of equation (\ref{epsilon}) and consequently of the proposition.
\end{proof}

An immediate consequence of Proposition \ref{same}  is the following corollary.

\begin{corollary}
For any graph $\Gamma$, $(Gr(\Gamma),\tau)$ is a graded, tracial, faithful $*$-probability space.\qed
\end{corollary}

We recall that if $(A,\tau)$ is a tracial probability space and $e \in A$
is a non-zero projection, then the corner $eAe$ is naturally a tracial probability space
where the trace is scaled so as to be $1$ on $e$. We will find
some corners of $Gr(\Gamma)$ to be useful. For a vertex
$v \in V$, we denote by $Gr(\Gamma,v)$ the probability space
$e_vGr(\Gamma)e_v$. Letting $e_0 = \sum_{v \in V_0} e_v$ - the
sum of the projections corresponding to the even vertices - we
will denote $e_0Gr(\Gamma)e_0$ by $Gr(\Gamma,0)$. 
Similarly, with $e_1 = \sum_{w \in V_1} e_w = 1-e_0$, we denote
$e_1Gr(\Gamma)e_1$ by $Gr(\Gamma,1)$.

The bipartite nature of the graph $\Gamma$ implies that the odd graded
pieces of these graded algebras reduce to zero. In particular,
$Gr(\Gamma,v)$ has as basis $\{[\xi] :$ $\xi$ is a path beginning and
ending at $v\}$ and is a connected graded algebra, while
$Gr(\Gamma,0)$ (resp. $Gr(\Gamma,1)$) has as basis $\{[\xi] :$ $\xi$
is a path beginning and ending at an even (resp. odd) vertex$\}$.
There are also corresponding notions in the $F(\Gamma)$ picture such
as $F(\Gamma,v)$, $F(\Gamma,0)$ or $F(\Gamma,1)$ and we will use
self-explanatory notation such as $P_n(\Gamma,0)$ or
$P_n(\Gamma,v)$. Thus, for instance, $F(\Gamma,v) = \oplus_{n \geq 0}
P_{2n}(\Gamma,v)$.
We will also tacitly use the fact that the isomorphism of $Gr(\Gamma)$ onto $F(\Gamma)$ of Proposition \ref{same} takes $Gr(\Gamma,v)$
to $F(\Gamma,v)$ for each $v \in V$. In particular, it takes $Gr(\Gamma,0)$
to $F(\Gamma,0)$ and $Gr(\Gamma,1)$
to $F(\Gamma,1)$.

Consider the Hilbert space $H(\Gamma)$ obtained by completing $F(\Gamma)$ for its inner-product, which has orthonormal
basis given by all 
\be\label{bracedef} \{\xi\} = \frac{1}{\sqrt{\mu(s(\xi))\mu(f(\xi))}} [\xi] \ee
where $\xi$ is a path in $\Gamma$. Equivalently, it is the Hilbert
space direct sum $\oplus_{n \geq 0} P_n(\Gamma)$ where
each $P_n(\Gamma)$ has orthonormal basis $\{\xi\}$
with $\xi$  a path in $\Gamma$ of length $n$.
We denote its norm by $|| \cdot ||$.

We also need the local Hilbert spaces $H(\Gamma,v)$ which we define to
be the completions of $F(\Gamma,v)$ for their trace norms. Note that
$F(\Gamma,v)$ is a (non-unital) subalgebra of $F(\Gamma)$ and that the
norm on $F(\Gamma,v)$ is a scaled version of the norm on $F(\Gamma)$.
The paths $[\xi]$ that begin and end at $v$ are an orthonormal
basis of  $F(\Gamma,v)$ (while they have norm $\mu(v)$ regarded
as elements of $F(\Gamma)$). 

We wish to show that the left regular representation of $F(\Gamma)$
on itself extends to a bounded representation on $H(\Gamma)$.
It clearly suffices to see that for $a \in P_m(\Gamma)$ and $b \in
F(\Gamma)$ there exists a constant $C$ (depending only on $a$)
such that $||a \# b|| \leq C||b||$. The proof of Proposition 4.3
in \cite{KdySnd2008} goes over to show that even the following is
sufficient (and that we may take $C = (2m+1)K$).

\begin{proposition}
For $a \in P_m(\Gamma)$ and $b \in
P_n(\Gamma)$ there exists a constant $K$ (depending only on $a$)
such that $||(a \# b)_t|| \leq K||b||$ for any $t$ with $|m-n| \leq t \leq m+n$.
\end{proposition}

\begin{proof} We will work with the orthonormal basis $\{\xi\}$ rather than
the orthogonal basis $[\xi]$. Observe that
$$
\{\xi\} \bullet \{\eta\} = 
\left\{ 
\begin{array}{ll}
                   0 & {\text {if\ }} f(\xi) \neq s(\eta) \\
                   \frac{1}{\mu(f(\xi))}{\{}\xi \circ \eta{\}} & {\text {if\ }} f(\xi) = s(\eta).
\end{array}
\right. 
$$
We may assume that $a = \{\xi\}$. 
Suppose that $b = \sum_{\eta} c_\eta \{\eta\}$, where the sum is over
all paths $\eta$ in $\Gamma$ of length $n$.
Since $(a \# b)_t$ is obtained by
an application of at most $m$ $S^n_i$'s
to $\sum_{\{\eta : f(\xi) = s(\eta)\}}  \frac{c_\eta}{\mu(f(\xi))} \{\xi \circ 
\eta\}$, it suffices to bound the operator norm of the $S^n_i$ and the 
Hilbert space norm of $\sum_{\{\eta : f(\xi) = s(\eta)\}}  \frac{c_\eta}{\mu(f(\xi))} \{\xi \circ 
\eta\}$.

It is easily checked that the adjoint of $S^n_i$ is given explicitly by
$$
(S^n_i)^*(\{\eta\}) = \sum_w \sum_{\rho : (v^\eta_{i-1} \stackrel{\rho}{\rightarrow} w) \in E} \frac{\mu(w)}{\mu(v^\eta_{i-1})} \{\eta_{[0,i-1]}\circ \rho \circ \tilde{\rho} \circ \eta_{[i-1,n-2]}\},
$$
and consequently that
$$
S^n_i(S^n_i)^*(\{\eta\}) = \left( \sum_w |\{{\rho : (v^\eta_{i-1} \stackrel{\rho}{\rightarrow} w) \in E\}}|   \left( \frac{\mu(w)}{\mu(v^\eta_{i-1})}\right)^2   \right)\{\eta\},
$$
for all $\{\eta\} \in P_{n-2}(\Gamma)$.
Thus, for a vertex $v$, if we define $$\delta(v) = \sum_w |\{{\rho : (v \stackrel{\rho}{\rightarrow} w) \in E\}}|    \left( \frac{\mu(w)}{\mu(v)}\right)^2,$$
and $\delta = max_{v \in V} \{\delta(v)\}$, then the operator
norm of $S^n_i$ is bounded above by $\sqrt{\delta}$.

Finally, note that $||b||^2 = \sum_\eta |c_\eta|^2$ while $$||\sum_\eta \frac{c_\eta}{\mu(f(\xi))} \{\xi \circ 
\eta\}||^2 \leq \frac{1}{\mu(f(\xi))^2} \sum_\eta |c_\eta|^2 =\frac{1}{\mu(f(\xi))^2}  ||b||^2.$$
Thus we may take $K = \frac{\max\{1,\delta^{m/2}\}}{\mu(f(\xi))}$ for
$a = \{\xi\} \in P_m(\Gamma)$ (the reason for the `max' being to allow
for the cases $\delta < 1$ and $\delta \geq 1$).
\end{proof}

We thus have a bounded left regular representation $\lambda :
F(\Gamma) \rightarrow {\mathcal L}(H(\Gamma))$ and we set $M(\Gamma) 
= \lambda(F(\Gamma))^{\prime\prime}$.
Similarly, for $i=0,1$, we have the left regular representation
$\lambda : F(\Gamma,i) \rightarrow {\mathcal L}(H(\Gamma,i))$ and we
may define $M(\Gamma,i) = \lambda(F(\Gamma,i))^{\prime\prime}$. 
It is easy to see that - see Lemma 4.4 of \cite{KdySnd2008} -
each of $M(\Gamma)$, $M(\Gamma,0)$ and $M(\Gamma,1)$ is a finite von Neumann
algebra.
The goal of the next section is to show that $M(\Gamma)$ is `almost
a $II_1$-factor'.

\section{Almost $II_1$-factoriality of $M(\Gamma)$}

Throughout this section, our standing assumption will be that the graph $\Gamma$ is connected and has at least one edge.
For such
a graph, it is clear that $Gr(\Gamma)$ is
infinite-dimensional.
The main results of this section apply only to graphs with at least two
edges and  show that  the
von Neumann algebra $M(\Gamma)$ is a direct sum of a $II_1$-factor and a  finite-dimensional abelian
algebra (possibly $\{0\}$) by analysing the local graded probability spaces
$Gr(\Gamma,v)$ for each vertex $v \in V$.

In the analysis of $Gr(\Gamma,v)$, an action by a certain category
that we denote by ${\mathcal C}(\delta)$ (where $\delta \in {\mathbb C}$ is some fixed non-zero parameter) will be extremely important, so we begin by describing this category.

Its objects are $[0],[1],[2],\cdots$, where we think of $[n]$ as a set
of $2n$ points on a horizontal line labelled $1,2,\cdots,2n.$ Note
that the objects of ${\mathcal C}(\delta)$ are denoted by exactly the
same notation as objects of ${\mathcal E}$ but mean different things.

The set $Hom([n],[m])$ is stipulated to have basis given by
all $T(P,Q)_{n}^{m}$ where $P \subseteq [m]$ and $Q \subseteq [n]$
are intervals of equal cardinality, where $T([4,5],[3,4])_5^8 \in
Hom([5],[8])$ is illustrated in Figure \ref{tab} below.
\begin{figure}[htbp]\label{tab}
\begin{center}
\includegraphics[height=2cm]{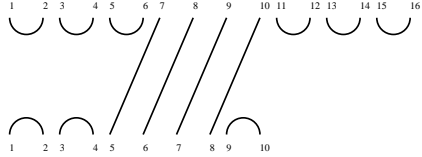}
\caption{The morphism $T([4,5],[3,4])_5^8 \in Hom([5],[8])$}
\end{center}
\end{figure}
The general prescription for $T(P,Q)_n^m$ is the following. Points below labelled by
elements of the sets $2Q$ and $2Q-1$ are joined to points above
labelled by $2P$ and $2P-1$ in order preserving fashion, and the rest
are capped or cupped  off without nesting.

Composition in ${\mathcal C}(\delta)$ for the basis elements is as in
Temperley-Lieb categories - by stacking the pictures and replacing
each closed loop that appears with a multiplicative $\delta$, thus
yielding a multiple of another basis element.  In order to explicitly
write down the composition rule in terms of the $T(P,Q)_n^m$, note
that there is a unique order preserving bijection $f_{PQ} : P
\rightarrow Q$ and that for $p \in P$, the marked points above that
are labelled by $2p-1$ and $2p$ are joined, respectively, to the
points below that are labelled by $2f_{PQ}(p)-1$ and $2f_{PQ}(p)$.  In the
example above, for instance, $f_{PQ}(p) = p-1$. Then the composition rule is
seen to be
\begin{equation*}
T(P,Q)^m_n \circ T(R,S)^n_p = \delta^{n  - |Q \cup R|} 
T(Y,Z)^m_p,
\end{equation*}
where $Y = f_{PQ}^{-1}(Q \cap R)$ and $Z = 
f_{RS}(Q \cap R)$.

For $n \geq 1$ let $A_-^n \in Hom([n],[n-1])$ be the morphism with
a single cap in the bottom left corner (i.e., $A_-^n = T(\{1,2,\cdots,n-1\},\{2,3,\cdots,n\})_n^{n-1}$) and $A_+^n \in Hom([n],[n-1])$ be the morphism with
a single cap in the bottom right corner (i.e., $A_+^n = T(\{1,2,\cdots,n-1\},\{1,2,\cdots,n-1\})_n^{n-1}$).
Similarly, for $n \geq 0$ let $C_-^n \in Hom([n],[n+1])$ be the morphism with
a single cup in the top left corner and $C_+^n \in Hom([n],[n+1])$ be the morphism with
a single cup in the top right corner. Thus, $C_-^n
=T(\{2,3,\cdots,n+1\},\{1,2,\cdots,n\})^{n+1}_n$ and $C_+^n
=T(\{1,2,\cdots,n\},\{1,2,\cdots,n\})^{n+1}_n$. (The letters $A$ and
$C$ are meant to suggest similarity to `annihilation' and creation'.)

\begin{proposition}\label{genrel}
If $\delta \neq 0$, the category ${\mathcal C}(\delta)$ is generated by the set
\[ \{ A^n_\pm: n \geq 1\} \cup \{ C^n_\pm: n \geq 0\}\]
of morphisms, and presented by the following relations, valid for all
$n \geq 0$:

\begin{eqnarray}
\label{smalla} A^1_- & = & A^1_+ \\
\label{smallc} C^0_- & = & C^0_+ \\
\label{amap} A^{n+1}_-  A^{n+2}_+  & = & A^{n+1}_+  A^{n+2}_- \\
\label{amcm} A^{n+1}_-  C^n_- & = & \delta~ id_{[n]} \\
\label{amcp} A^{n+2}_-C^{n+1}_+ & = & C^n_+A^{n+1}_- \\
\label{apcm} A^{n+2}_+C^{n+1}_- & = & C^n_-A^{n+1}_+ \\
\label{apcp} A^{n+1}_+  C^n_+ & = & \delta ~id_{[n]}\\
\label{cmcp} C^{n+1}_-  C^n_+  & = & C^{n+1}_+  C^n_- 
\end{eqnarray}

More explicitly, suppose a category ${\mathcal D}$ has the property
that $Hom(D, D^\prime)$ is a complex vector space for every pair
$(D,D^\prime)$ of 
  objects. Then, in order to specify a functor from ${\mathcal
    C}(\delta)$ to ${\mathcal D}$, it (is necessary and it) suffices
  to find objects $D_n \in {\mathcal D}$ for $n \geq 0$ and morphisms
  $\widetilde{A}^n_\pm : D_n \rightarrow D_{n-1}$ for $n \geq 1$
  and $\widetilde{C}^n_\pm
  : D_n \rightarrow D_{n+1}, $ for $n \geq 0$ satisfying the
  relations (\ref{smalla})-(\ref{cmcp}) above.\qed

\end{proposition}

\begin{proof}
It is easy to check that the relations (\ref{smalla})-(\ref{cmcp}) are satisfied.

To see that these are the only relations, we need to first observe that
the following identities, for $k,l \geq 0$, are consequences of them:

\begin{eqnarray} 
A^1_- \cdots A^k_-A^{k+1}_+ \cdots A^{k+l}_+ &=& A^1_+ \cdots
A^{k+l}_+ \label{onlya} \\
C^{k+l-1}_+ \cdots C^k_+ C^{k-1}_- \cdots C^0_- &=& C^{k+l-1}_+ \cdots C^0_+ \label{onlyc}
\end{eqnarray}

These two identities are seen inductively to follow from the equations
numbered (\ref{smalla})  and (\ref{amap}), and from (\ref{smallc}) and (\ref{cmcp}) respectively, from among
the above relations.

We next describe a `canonical form' for every morphism in ${\mathcal C}(\delta)$
as a word in the generators, in such a way that if we assign
the `rank'  1, 2, 3, and 4 to any $C_+, C_-, A_+$ and $A_-$ respectively,
then if the word contains generators of ranks $i$ and $j$, with $i<j$,
then the generator of rank $i$ will appear to the left of the one with
rank $j$. For example, the morphism illustrated in Figure \ref{tab})
is expressed as 
\[T([4,5],[3,4])_5^8 = C^7_+C^6_+C^5_+C^4_-C^3_-C^2_-A^3_+A^4_+A^5_-\]
(The algorithm for arriving at this word is to `first list all the
caps from left to right, and then all the cups also from left to
right'.)  Notice that it is only when there are no through strings
that there is an ambiguity (about whether to choose a + or a - for the
$A$'s and $C$'s, but this is resolved using equations (\ref{onlya})
and (\ref{onlyc}) above, using which we can demand in the case of no
through strings, that all $C$'s and $A$'s come with the subscript `+'. On
the other hand, if there are through strings, the number of through
strings can be read off from this word, at the point of transition
from $C$'s to $A$'s (the number of through strings is exactly twice
the superscript of the rightmost $C$). (By the way, it is in order to
lay hands on $id_{[n]}$ that we need the condition $\delta \neq 0$.)

Finally, if the `rank ordering' specified above is violated in any
word in the generators, such instances may be rectified uniquely by using
(\ref{amap})-(\ref{cmcp}).
(For example, any instance where an $A_-$ (of rank 4) precedes any
generator of rank 3, 2, or 1, is set right by equations numbered (\ref{amap})-(\ref{amcp}).)
\end{proof}

We use Proposition \ref{genrel} to get a functor from ${\mathcal
  C}(\delta(v))$ to the category ${\mathcal D}$ of $\C$-vector spaces
and $\C$-linear maps.
Recall that $\delta(v)  = \sum_w |\{{\rho : (v \stackrel{\rho}{\rightarrow} w) \in E\}}|    \left( \frac{\mu(w)}{\mu(v)}\right)^2.$
Let
$D_n$ be $P_{2n}(\Gamma,v)$ and define the action of the generating morphisms
as follows on $[\xi] \in P_{2n}(\Gamma,v)$.
\begin{eqnarray*}\label{action}
\begin{array}{ll}
\widetilde{A}^n_- ([\xi]) = \delta_{\xi_1,\widetilde{\xi_2}} \frac{\mu(v^\xi_1)}{\mu(v)} \xi_{[2,2n]}    & 
\widetilde{C}^n_- ([\xi]) = \sum_w \sum_{\rho : (v \stackrel{\rho}{\rightarrow} w) \in E} \frac{\mu(w)}{\mu(v)} [ \rho \circ \tilde{\rho} \circ \xi] \\
\widetilde{A}^n_+ ([\xi]) = \delta_{\xi_{2n-1},\widetilde{\xi_{2n}}} \frac{\mu(v^\xi_{2n-1})}{\mu(v)} \xi_{[0,2n-2]}  & 
\widetilde{C}^n_+ ([\xi]) = \sum_w \sum_{\rho : (v \stackrel{\rho}{\rightarrow} w) \in E} \frac{\mu(w)}{\mu(v)} [\xi \circ  \rho \circ \tilde{\rho}]
\end{array}
\end{eqnarray*}

A little calculation now proves the following.
\begin{proposition}
The action by the generators given by the equations above extends
to give a well defined functor from ${\mathcal C}(\delta(v))$ to 
${\mathcal D}$.\qed
\end{proposition}

Note that regarding $P_{2n}(\Gamma,v)$ as subspaces of
$H(\Gamma,v)$, the maps $\widetilde{A}^{n+1}_-$ and $\widetilde{C}^{n}_-$ are adjoints of each other, as are
$\widetilde{A}^{n+1}_+$ and $\widetilde{C}^{n}_+.$

We now try to determine the structure of the center of $Gr(\Gamma,v)$.
In particular, we show that it is at most two dimensional. We will find the following notation useful. For $a \in F(\Gamma,v)$ let
$[a] = \{ \xi \in H(\Gamma,v) : \lambda(a)(\xi) = \rho(a)(\xi)
\},$\footnote{As usual, we wreite $\rho(a) = J\lambda(a)^*J$, with $J$
  being the modular conjugation operator.}
which is a closed subspace of $H(\Gamma,v)$.
Denote by $\Omega$ the vector $e_v \in F(\Gamma,v) \subseteq 
H(\Gamma,v)$, note that this is (cyclic and) separating for $M(\Gamma,v)$
and thus the operator equation $ax = xa$ is equivalent to the
vector inclusion $x\Omega \in [a]$ for $x \in M(\Gamma,v)$.

Our strategy is similar to that in \cite{KdySnd2008}, with some  differences. 
We first define two elements
$c \in P_2(\Gamma,v)$ and $d \in P_4(\Gamma,v)$. 
(For notational convenience we do not use the possibly more correct
notation $c_v$ and $d_v$.) We then
show that $[c] \cap [d]$ is 1-dimensional if $\delta(v) \geq 1$
and 2-dimensional if $\delta(v) <1$ (assuming that the graph $\Gamma$ in
question has at least two edges). Finally we show that in case $\delta(v) < 1$, the centre
of $Gr(\Gamma,v)$ is actually 2-dimensional and that the cut-down
by one of the central projections is just ${\mathbb C}$. Now some
simple computations give the desired result that $M(\Gamma)$ is
either a $II_1$-factor or a direct sum of one with a finite-dimensional
abelian algebra.

In the sequel, we dispense with `tilde's and continue to use the same
symbol for morphisms in ${\mathcal C}(\delta)$ and the associated linear maps
between the $P_{2n}(\Gamma,v)$.

Let $c \in P_2(\Gamma,v)$ be the element $C^0_-(1)$. Explicitly,
$$c = \sum_w \sum_{\rho : (v \stackrel{\rho}{\rightarrow} w) \in E} \frac{\mu(w)}{\mu(v)} [\rho \circ \tilde{\rho}].$$ 
Let $c_0 =1$ and by 
$c_{2n}$ for $n \geq 1$, we will denote the
element $C^{n-1}_-C^{n-2}_-\cdots C^0_-(1).$ Thus $c_{2n} \in
P_{2n}(\Gamma,v)$ and $c_2 = c$. Note that by induction on $n$, $c_{2n}$
is seen to be the highest degree 
term of $c^n$ in $F(\Gamma,v)$ and  to be a polynomial
in $c$ of degree $n$.
Let $C \subseteq H(\Gamma,v)$
be the closed subspace spanned by all the $c_{2n}\Omega$ for $n \geq
0$. We then have the following crucial result 
which is the analogue of Proposition 5.4 of \cite{KdySnd2008}.

\begin{proposition}\label{c=C}
$[c] = C$.
\end{proposition}

Before sketching a proof, we state a key lemma used which is the
analogue of Lemma 5.6 of \cite{KdySnd2008}. By $C_{2n}$ we denote the
(1-dimensional) subspace of $P_{2n}(\Gamma,v)$ spanned by 
$c_{2n}$ and by $C_{2n}^\perp$ its orthogonal complement in $P_{2n}(\Gamma,v)$. Thus if $C^\perp$ is the orthogonal complement of
$C$ in $H(\Gamma,v)$, then $C^\perp = \oplus_{n \geq 0} C_{2n}^\perp$.

\begin{lemma}\label{ccommlem}
For $n \geq 0$, the map $C_{2n}^\perp \ni x \mapsto z = (c\#x - x \#
c)_{2n+2} \in P_{2n+2}(\Gamma,v)$ 
is injective with inverse
given by 
$$
x = \sum_{t=1}^{n} \delta(v)^{-t} T([1,n+1-t],[t+1,n+1])^n_{n+1}(z).
$$
\end{lemma}

We omit the proof except to remark that $(c\#x)_{2n+2}$ and
$(x\#c)_{2n+2}$ are just $C^n_-(x)$ and $C^n_+(x)$. We also omit the
proof of the next corollary which is the analogue of Corollary 5.7 of
\cite{KdySnd2008} with identical proof. 

\begin{corollary}\label{xnxm}
Suppose that $\xi = (x_0,x_{1},\cdots) \in \oplus_{n=0}^\infty C_{2n}^{\perp} =  C^{\perp}$ and satisfies $\lambda(c)(\xi) = \rho(c)(\xi)$.
i.e., $\xi \in C^\perp \cap [c]$. Then, for $m > n > 0$ with $m-n = 2r$, we have:
\begin{eqnarray*}
x_n &=& \sum_{t=1}^{n} \delta(v)^{-(t+r-1)}  \left\{T([1,n+1-t],[t+r,n+r])^{n}_{m}(x_m)\right.\\
 & &  - \left.T([1,n+1-t],[t+r+1,n+r+1])^{n}_{m}(x_m)\right\}
\end{eqnarray*}
\end{corollary}

One more result needed in proving Proposition \ref{c=C} is the following norm estimate which is the analogue of Lemma 5.8 of
\cite{KdySnd2008}.
\begin{lemma}\label{easy}
Suppose that $x = x_n \in P_{2n}(\Gamma,v) \subseteq F(\Gamma,v)$ and let $y = T(P,Q)^m_n(x) \in P_{2m} \subseteq F(\Gamma,v)$ for some morphism
$T(P,Q)^m_n$. Then $||y|| \leq \delta(v)^{\frac{1}{2}(n+m) - |P|} ||x||$.
\end{lemma}

\begin{proof}
Consider the linear extension of the map defined by $Hom([n],[m]) \ni
T(P,Q)^m_n \mapsto  \delta(v)^{\frac{1}{2}(n+m) - |P|} \in \C$ for
each $n,m \geq 0$. 
The observation is that it is multiplicative on composition. For consider
$T(P,Q)^m_n \in Hom([n],[m])$ and $T(R,S)^n_p \in Hom([n],[p])$.
Their composition is given by $\delta(v)^{n  - |Q \cup R|} 
T(Y,Z)^m_p$ where $Y = f_{PQ}^{-1}(Q \cap R)$ and $Z = 
f_{RS}(Q \cap R)$.
The multiplicativity assertion amounts to verifying that
$$
\frac{1}{2}(n+m) - |P| + \frac{1}{2}(n+p) - |R| = n  - |Q \cup R| +  \frac{1}{2}(m+p) - |Y|,
$$
which is easily verified to hold.

Hence it suffices to verify the norm estimate when $T(P,Q)^m_n$
is one of the generators $C^n_\pm$ or $A^n_\pm$ of the category ${\mathcal C}(\delta(v))$.
Note that the norm estimate for all these generators is just $||y|| \leq \delta(v)^{\frac{1}{2}} ||x||$.
For $C^n_\pm$, we have $(C^n_\pm)^*C^n_\pm = A^{n+1}_\pm C^n_\pm = \delta(v) I^n$ while for $A^n_\pm$, we have $A^n_\pm (A^n_\pm)^* = A^n_\pm C^{n-1}_\pm = \delta(v) I^{n-1}$.
This proves the norm estimate for generators and completes the proof of the lemma.
\end{proof}

\begin{proof}[Proof of Proposition \ref{c=C}]
Each $c_{2n}$, being a polynomial in $c$, clearly commutes with $c$
and it follows that $C \subseteq [c]$. To prove the other inclusion, it is enough to see that $C^\perp \cap [c] = \{0\}$. Suppose that
$\xi = (x_0,x_{1},\cdots) \in C^\perp \cap [c]$. Note that $x_0 =0$
since $x_0 \in C^\perp_0$ while $C_0 = P_0(\Gamma,0) ={\mathbb C}$.
By Corollary \ref{xnxm}, we have for $m > n >0$ with $m-n = 2r$,
\begin{eqnarray*}
x_n &=& \sum_{t=1}^{n} \delta(v)^{-(t+r-1)}  \left\{T([1,n+1-t],[t+r,n+r])^{n}_{m}(x_m)\right.\\
 & &  - \left.T([1,n+1-t],[t+r+1,n+r+1])^{n}_{m}(x_m)\right\}
\end{eqnarray*}
Now, applying the norm estimate from Lemma \ref{easy} and using the
triangle inequality gives $||x_n|| \leq 2n||x_m||$. Now $\xi \in
H(\Gamma,v) \Rightarrow lim_{m \rightarrow \infty} ||x_m|| = 0$ and so $x_n =0$
for all $n >0$. Hence $\xi = 0$, completing the proof.
\end{proof}

\begin{figure}[htbp]
\begin{center}
\includegraphics[width=1cm]{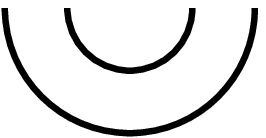}
\caption{Pictorial representation of $d$}\label{Dfig}
\end{center}
\end{figure}

We now consider the element $d \in P_4(\Gamma,v)$ defined explicitly by
$$
d = 
\sum_w
\sum_{\rho : (v \stackrel{\rho}{\rightarrow} w) \in E}
\sum_x
\sum_{\zeta : (w \stackrel{\zeta}{\rightarrow} x) \in E}
\frac{\mu(x)}{\mu(v)}
[
v \stackrel{\rho}{\rightarrow} w \stackrel{\zeta}{\rightarrow} x 
\stackrel{\widetilde{\zeta}}{\rightarrow} w \stackrel{\widetilde{\rho}}{\rightarrow} v
].
$$
Loosely, $d$ can be thought of as pictorially represented as in Figure \ref{Dfig}.

It must be noted that the action of the category
${\mathcal C}(\delta)$ on $F(\Gamma,v)$ may be extended to an action of the
entire Temperley-Lieb category on $F(\Gamma,v)$ if
$\mu^2$ is a Perron-Frobenius eigenvector for
the incidence matrix of $\Gamma$, however, this may no longer be
possible if the `Perron-Frobenius assumption' is dropped. Since we
shall have to address that 
situation, we do not assume this. Nevertheless, we will see that in
situations of interest 
for us, this pictorial representation will be of heuristic value.

We wish to consider a special
element of $M(\Gamma,v)$ defined when $\delta(v) < 1$.  Since the
proof of existence of this element requires a careful norm estimate,
we digress with the necessary lemma.

\begin{lemma}\label{zdef}
  Suppose the weighting on  $\Gamma$ is such that
  $\delta(v) < 1$.  The sequence of elements $x_m = \sum_{n=0}^m
  (-1)^n c_{2n} \in F(\Gamma,v) \subseteq \footnote{This inclusion is,
  naturally, via $\lambda$.}
 M(\Gamma,v) \subseteq
  {\mathcal L}(H(\Gamma,v))$ converge in the strong operator
  topology. Hence the series $\sum_{n=0}^\infty (-1)^n c_{2n}$ defines an
  element $z_v \in M(\Gamma,v).$
\end{lemma}

\begin{proof}
  It suffices to see that the $x_m$ are uniformly bounded in norm and
  that for $\xi$ in a dense subspace of $H(\Gamma,v)$, the sequence
  $x_m\xi$ converges in $H(\Gamma,v)$. Note that if $\xi = \Omega
  (~=[(v)]~)$ - the vacuum vector of $H(\Gamma,v)$ - then, $x_m\xi =
  \sum_{n=0}^m (-1)^n c_{2n}\Omega.$ Since $||c_{2n}\Omega||^2 =
  \delta(v)^n$, when $\delta(v) < 1$, the $x_m\xi$ converge in
  $H(\Gamma,v)$. It follows that on the dense subspace
  $M(\Gamma,v)^\prime \Omega$ too, we obtain convergence. It remains
  to prove the norm estimate.

Consider the block matrix  representing left-multiplication
by $c_{2n}$ on the Hilbert space $H(\Gamma,v)$ with respect
to the orthogonal decomposition $H(\Gamma,v) = \oplus_{n=0}^\infty
P_{2n}(\Gamma,v)$. The definition of multiplication in $F(\Gamma,v)$
shows that for any path $\xi$ of length $2j$,
$$
c_{2n} \# [\xi] = \sum_{k=0}^{min\{2j,2n\}} C^{n-\lceil \frac{k}{2} \rceil}A^{\lfloor \frac{k}{2} \rfloor} [\xi]
$$
where, to avoid heavy notation, for $[\xi] \in P_{2j}(\Gamma,v)$ we write $C^pA^q [\xi]$  to mean
$$
C_-^{j-q+p-1}\cdots C_-^{j-q+1}C_-^{j-q}A^{j-q+1}_- \cdots A^{j-1}_-A^j_- [\xi].
$$
Thus the $(i,j)$-block (note that $i,j \geq 0$) of the matrix for
$\lambda(c_{2n})$  
is $0$ unless $|i-j| \leq n \leq i+j$, in which case it is given by
$$
C^{\lfloor \frac{n-j+i}{2} \rfloor} A^{\lfloor \frac{n+j-i}{2} \rfloor}.
$$
It now follows that the  $(i,j)$-block of the matrix for $x_m$ is  given
by
$$
\sum_{|i-j| \leq n \leq min\{m,i+j\}} (-1)^n C^{\lfloor \frac{n-j+i}{2} \rfloor} A^{\lfloor \frac{n+j-i}{2} \rfloor}
$$

Every odd term of this sum (starting with the first which corresponds to
$n = |i-j|$) is equal, except for sign, with the succeeding even term and
so the sum vanishes when there are an even number of terms and
equals its last term when there are an odd number of terms. 
Note that the number of terms is $min\{m,i+j\} - |i-j|+1$.

We consider two cases depending on whether $m \leq i+j$ or $m > i+j$.\\
Case I: If $m > i+j$, the number of terms is certainly odd and so
$$
(x_m)^i_j = (-1)^{i+j} C^{i}A^j
$$
Case II: If $i+j \geq m$ the number of terms is odd or even according
as $m$ and $i+j$ have the same parity or not and so, in this case,
$$
(x_m)^i_j = 
\left\{
\begin{array}{ll}
0 &{\text {if $m < |i-j|$ or $(i+j) -m$ is odd}} \\
 (-1)^m C^{\lfloor \frac{m-j+i}{2} \rfloor} A^{\lfloor \frac{m+j-i}{2} \rfloor} &{\text {otherwise}}
\end{array}
\right.
$$

For instance, the matrix for $x_2$ is given by
$$
\left(
\begin{array}{rrrrrr}
I & -A & A^2 & 0 & 0 & \cdots\\
-C &CA & 0 & A^2 & 0 & \cdots\\
C^2 & 0 & CA & 0 & A^2 & \cdots\\
0 & C^2 & 0 & CA & 0 & \cdots\\
0 & 0 & C^2 & 0 & CA &  \cdots\\
\vdots & \vdots & \vdots & \vdots & \vdots &  \ddots
\end{array}
\right)
$$

Observe that the $(i,j)$ entry of $x_m$ is non-zero only if $|i-j| \leq m$
in which case it is of the form $\pm C^pA^q$ with $p+q \geq |i-j|$.
Since each of $A$ and $C$ has norm bounded above by $\delta(v)^{\frac{1}{2}}$, the diagonal of $x_m$ with $i-j =t$ has norm at most
$\delta(v)^\frac{|t|}{2}$ and so $x_m$ itself has norm bounded by $\sum_{t=-m}^m \delta(v)^\frac{|t|}{2} \leq 1 + 2 \sum_{t=1}^\infty \delta(v)^\frac{t}{2}$. Thus the $x_m$ are uniformly bounded in norm, finishing the proof.
\end{proof}

For reasons that will  become clear in Proposition \ref{1or2}, we define (in the notation
introduced in Lemma \ref{zdef}),
\[ e_v^I = \left\{ \ba{ll} 0 & \mbox{if } \delta(v) \geq 1\\
(1-\delta(v)) z_v & \mbox{if } \delta(v) < 1 \ea \right. \]
and $e_v^{II} = e_v - e_v^I$.

The next proposition is the analogue of
Proposition 5.5 of \cite{KdySnd2008}; the reader is urged to
compare the proof of that proposition with this one.
For the rest of this section,  $\Gamma$ will always denote a connected graph with at least two edges.

\begin{proposition}\label{ccapd}
  Suppose that $\Gamma$ is a connected graph with at least two edges.
 Then, $[c] \cap [d]$ has basis $\{1\}$ if $\delta(v) \geq 1$
  and basis $\{1,z_v\}$ if $\delta(v) < 1$.
\end{proposition}

\begin{proof} Suppose that $\xi = (x_0,x_2,x_4,\cdots) \in [c] \cap [d]$.
Since $\xi \in [c] = C$ by Proposition~\ref{c=C}, there exist scalars $y(2n) \in {\mathbb C}$
such that $x_{2n} = y(2n)c_{2n}$ and since $\xi \in H(\Gamma,v)$, we have that $||\xi||^2= \sum_{n=0}^\infty |y(2n)|^2 \delta(v)^n < \infty$.

Some computation now shows that for $t>2$, the $P_{2t}(\Gamma,v)$
component of $\lambda(d)([\xi]) - \rho(d)([\xi])$ is given by
$(y(2t-2) + y(2t-4))(d\#c_{2t-4} - c_{2t-4}\#d)_{2t}$. Writing out
$(d\#c_{2t-4} - c_{2t-4}\#d)_{2t}$ for $t>2$ in terms of the path basis, inspection
shows that it cannot vanish since $\Gamma$ has at least two edges.
Thus each $y(2t-2) + y(2t-4)$ vanishes for $t >2$.

Hence if $y(2) = y$, then all $y(4n-2) =y$ and all $y(4n) = -y$
(for $n\geq 1$). Now the norm condition shows that if $\delta(v) \geq 1$, then $y=0$
so that
$\xi$ is a multiple of $(1,0,0,\cdots)$, while
if $\delta(v) < 1$, then $\xi$ is  a linear combination of $1$ and
$z_v$.
\end{proof}

We now justify the choice of notation for $e_v^I$ and $e_v^{II}$.

\begin{proposition}\label{1or2}
Suppose that $\Gamma$ is a connected graph with at least two edges. Then,
\ben
\item If $\delta(v) < 1$,  then $Z(M(\Gamma,v))$ is 2-dimensional
and has basis $\{e^I_v, e^{II}_v\}$. The element $e_v^I \in
M(\Gamma,v)$ is a minimal (and central) projection.
\item $e^{II}_v$ is a minimal central projection in $M(\Gamma,v)$ and
  $e^{II}_vM(\Gamma,v)e^{II}_v$ is a $II_1$-factor; in particular,
  $M(\Gamma,v)$ is a factor if $\delta(v) \geq 1$.
\een
\end{proposition}

\begin{proof} (1) Since $Z(M(\Gamma,v))\Omega \subseteq [c] \cap [d]$ and
  $Z(M(\Gamma,v)) \rightarrow Z(M(\Gamma,v))\Omega$ is injective,
  Proposition \ref{ccapd} implies that $Z(M(\Gamma,v))$ is at most
  2-dimensional.

  If $\delta(v) <  1$, it follows from Proposition \ref{ccapd} that
  $Z(M(\Gamma,v))$ is at most two dimensional and contained in the
  span of $e_v$ (the identity for $M(\Gamma,v)$) and $e_v^I$.

We shall regard $c_{2n}$ as an element of $e_v F(\Gamma) e_v
  \subset F(\Gamma)$ (so that $c_0 = e_v$) and consider a length one path
  $\xi$ of $F(\Gamma)$.
  Calculation shows that
$$
c_{2n} \# [\xi] = 
c_{2n} \bullet [\xi] + c_{2n-2} \bullet [\xi] 
$$
and therefore that
$$
z_v \# [\xi] = 0
$$
Similarly, (or by simply taking adjoints) we find that also
$$
[\xi] \# z_v  = 0
$$
Associativity of multiplication implies that these equations hold even
when $[\xi]$ is a path of length greater than $1$ from $w$ to $x$.

Finally, if $\xi$ is a path of length 0, then $\xi = e_w$ for some
$w$, and
\[ z_v \# [\xi] = [\xi] \# z_v  = \left\{ \ba{ll} 0 & \mbox{if } w \neq v\\z_v & \mbox{if } w
  = v \end{array} \right. ~.\]
It follows easily that $z_v \in Z(M(\Gamma))$ and that 
\bea \label{e^I_v}z_v M(\Gamma) = \C z_v
\eea

Taking adjoints yields 
\bea \label{e^I_v2}M(\Gamma) z_v = \C z_v
\eea
Deduce that there is some constant
$\gamma > 0$ such that $z_v^2 = \gamma z_v$.

By comparing their inner products (in $\CH(\Gamma,v)$) with $e_v$, we
find that
\[ \sum_{n=0}^\infty (\delta(v))^n = \gamma\]
and hence that $\gamma = (1 - \delta(v))^{-1}$. Thus we find that
indeed $e^I_v$ is a projection in $M(\Gamma)$ which is central and
minimal, since 
\be \label{eImin} e_v^I M(\Gamma) e_v^I = \C e_v^I. \ee

\medskip \noindent (2) 
It is seen from
Proposition \ref{ccapd} that in both cases, $e_v^{II}$ is a minimal
central projection. The assumed non-triviality of $\Gamma$ ensures
that $M(\Gamma,v)$ is an infinite-dimensional 
but finite von Neumann algebra, it is seen (from the minimality of
$e^I_v$ in case $\delta(v) < 1$) that the localisation
$e^{II}M(\Gamma,v)e^{II}$ is a $II_1$ factor.
\end{proof}

\begin{corollary}\label{Isorth}
For distinct vertices $v$ and $w$, we have
\[e^I_v M(\Gamma) e_w = \{0\} = e_w M(\Gamma) e^I_v\]
\end{corollary}

\begin{proof} We only need to prove the nontrivial case when
  $\delta(v) < 1$.
First deduce from equation (\ref{e^I_v}) (when $\delta(v) < 1$) and
the definition of $e^I_v$ 
(when $\delta(v) \geq 1$) that $e^I_v M(\Gamma) = \C e^I_v$. Hence,
  $e^I_v M(\Gamma) e_w = \C e^I_v e_w = 0$. The second assertion of the
  corollary is obtained by taking adjoints in the first.
\end{proof}

Before proceeding to the next corollary to Proposition \ref{1or2} we digress with an elementary fact about local and global behaviour of von Neumann algebras.

\begin{lemma}\label{locglob}
Suppose $\{p_i : i \in I\}$ is a partition of the identity element  into a family of pairwise orthogonal projections of a von Neumann algebra $M$.

Then the following conditions are equivalent:
\ben
\item $M$ is a factor;
\item $p_iMp_i$ is a factor for all $i \in I$, and $p_iMp_j \neq \{0\} ~\forall i,j~.$
\een
\end{lemma}

\begin{proof}
We only indicate the proof of the non-trivial implication $(2) \Rar (1)$.
For this, suppose $x \in Z(M)$. Let us write $x_{ij} = p_ixp_j$. The
assumption (2) clearly implies that (i) $x_{ij} = 0$ for $i \neq j$
(since $x$ commutes with each $p_i$); and (ii) for each $i \in I$,
$x_{ii} = \lambda_i p_i$ for some $\lambda_i \in \C$. Fix an arbitrary
pair $(i,j)$ of distinct indices from $I$. By assumption, there exists
a non-zero $y \in M$ satisfying $y = p_iyp_j$. The requirement $xy =
yx$ is seen to now imply that $\lambda_i y = \lambda_j y$ and hence
that  $\lambda_i = \lambda_j$, and this is true for all $i,j$. Hence,
$ x = \sum_i \lambda_i p_i \in \C 1$.
\end{proof}

\begin{corollary}\label{IIseq}
Assume that $\Gamma$ is a connected graph with at least two edges. Then,
\ben
\item $e^{II}_vM(\Gamma)e^{II}_w \neq 0$ for all vertices $v,w$
\item If we let $e^{I} = \sum_{v \in V} e^{I}_v$, then
$e^{I}M(\Gamma)e^{I} = \bigoplus_{v \in V} \C e^I_v$.
\item If we let $e^{II} = \sum_{v \in V} e^{II}_v$, then
  $e^{II}M(\Gamma)e^{II}$ is a $II_1$ factor.
\een
\end{corollary}

\begin{proof} 
(1) Since $\Gamma$ is connected, we can find a finite path $\xi$ with
$s(\xi) = v$ and $f(\xi) = w$. Then, deduce from Corollary
\ref{Isorth} that
\beast
0 &\neq& [\xi]\\
 &=& e_v \# [\xi] \# e_w\\
&=& e^{II}_v \# [\xi] \# e^{II}_w
\eeast

(2) First observe that as $e^I_w \leq e_w$ for all $w \in V$, it
follows that
\beast
v \neq w & \Rar & e^I_vM(\Gamma)e^I_w\\
&=& e^I_vM(\Gamma) e_w e^I_w\\
&=& 0
\eeast
by Corollary \ref{Isorth}.

On the other hand
\[ e^I_v M(\Gamma)e^I_v = \C e^I_v\]

The desired assertion follows from the orthogonality of the $e^I_v$'s.

(3) As has already been observed, $e_v M(\Gamma) e_v = M(\Gamma,v)$
since $M(\Gamma)$ (resp., $M(\Gamma,v)$) is generated (as a von
Neumann algebra) by the set of all finite paths $[\xi]$ (resp., those
paths which start and finish at $v$, i.e., which satisfy $[\xi] = e_v
\# [\xi] \# e_v)$. Given this observation, the assertion to be proved
is seen to follow from Proposition \ref{1or2}, Lemma \ref{locglob} and
the already established part (1) of this Corollary. 
\end{proof}

We have finally arrived at the main result of this section, whose
statement uses the foregoing notation and which is an immediate
consequence of Corollary \ref{IIseq}.

\begin{theorem}\label{Mgamstr}
Assume $\Gamma$ is a connected graph with at least two edges.
Then, we have the following isomorphism of non-commutative
probability spaces:
\[M(\Gamma) \cong M \oplus \bigoplus_{\{v \in V: \delta(v) < 1\}} \underset{(1-\delta(v)) \mu^2(v)}{\C  e^I_v}\]
where $M$ is some $II_1$ factor (and we have omitted mention of the
obvious value of the trace-vector on the $M$-summand for typographical reasons).
\end{theorem}

\begin{proof}
It follows from Proposition \ref{1or2}(1), Corollary \ref{Isorth} and the fact that $1 =
\sum_v e_v = \sum_v e^I_v + e^{II}$ that each $e^I_v \in Z(M(\Gamma))$ and that consequently both
$e^I$ and $e^{II}$ are central. The asserted conclusion follows now from Corollary \ref{IIseq}.
\end{proof}

\begin{corollary}\label{mgamfac}
If $~ \Gamma$ is as in Theorem \ref{Mgamstr}, and is equipped with the
`Perron-Frobenius weighting', then $M(\Gamma)$ is a $II_1$ factor.
\end{corollary}

\begin{proof}
The hypotheses ensure that for all $v$, $\delta(v) = \delta$ is the
Perron-Frobenius eigenvalue of the adjacency matrix of $\Gamma$, which
in turn is greater than one, so the second summand of Theorem
\ref{Mgamstr} is absent.
\end{proof}

\section{Structure of the even graded probability space}

In this section we let $\Gamma$ be any finite, weighted, bipartite graph and regard $Gr(\Gamma,0)$ as an operator valued probability space
over its subalgebra $P_0(\Gamma,0)$ - the abelian algebra with
minimal central projections given by all $e_v$ where $v \in V_0$
is an even vertex. Our goal is to express this as a(n algebraic) free product
with amalgamation over $P_0(\Gamma,0)$ of simpler subalgebras.

We briefly summarise from \cite{Spc1998} the theory of operator valued probability spaces and  operator valued free cumulants.
An operator valued probability space is a unital inclusion of unital algebras $B \subseteq A$
equipped with a $B-B$-bimodule map $\phi : A \rightarrow B$ with $\phi(1) =
1$. 
A typical example is $N \subseteq M$ where $M$ is a von Neumann
algebra with a faithful, normal, tracial state $\tau$ and
$\phi$ is the $\tau$-preserving conditional expectation.

The lattice of non-crossing partitions plays a fundamental role
in the definition of free cumulants. Recall that for a totally ordered finite set $\Sigma$,
a partition $\pi$ of $\Sigma$ is said to be non-crossing if
whenever $i<j$ belong to a class of $\pi$ and $k<l$
belong to a different class of $\pi$, then it is not the case that
$k < i < l < j$ or $i < k < j < l$. The collection of non-crossing partitions of $\Sigma$, denoted
$NC(\Sigma)$, forms a lattice for the partial order defined by $\pi \geq  \rho$
if $\pi$ is coarser than $\rho$ or equivalently, if $\rho$ refines $\pi$.
The largest element of the lattice $NC(\Sigma)$  is denoted
$1_\Sigma$. Explicitly, $1_\Sigma = \{\Sigma\}$. If $\Sigma = [n] \stackrel{\mbox {\tiny {def}}}{=} \{1,2,\cdots,n\}$
for some $n \in {\mathbb N}$,  we will write
$NC(n)$ and $1_n$ for $NC(\Sigma)$ and $1_\Sigma$ respectively.

Before defining operator valued free cumulants, we state a basic
combinatorial result  that we will refer to as M\"{o}bius inversion.
Suppose that $A$ is an operator valued probability space over $B$.
Let $X \subseteq A$ be a $B-B$-submodule and suppose given $B-B$-bimodule maps
$\phi_n : \otimes_B^n X \rightarrow B$. By the multiplicative
extension of this collection, we will mean the collection of
$B-B$-bimodule maps $\{\phi_\pi : \otimes_B^n X \rightarrow B\}_{n
  \in {\mathbb N}, \pi \in NC(n)}$ defined recursively by 

\begin{eqnarray*}
\lefteqn{\phi_\pi(x^1\otimes x^2\otimes \cdots\otimes x^n) = }\\
&\left\{
\begin{array}{ll}
\phi_n(x^1\otimes x^2\otimes \cdots\otimes x^n) &\\
\phi_\rho(x^1\otimes \cdots\otimes x^{k-1}\otimes x^{k}\phi_{l-k}(x^{k+1}\otimes \cdots\otimes x^l)\otimes x^{l+1}\otimes \cdots\otimes x^n), 
\end{array}
\right.
\end{eqnarray*}
according as $\pi=1_n$ or $\pi = \rho \cup 1_{[k+1,l]}$ for $\rho \in NC([1,n] \setminus [k+1,l])$.
A little thought shows that the multiplicative extension is well-defined.
Let $\mu(\cdot,\cdot)$ be the M\"{o}bius function of the lattice $NC(n)$ - see Lecture 10 of \cite{NcaSpc2006}.

\begin{proposition}\label {Mobius}
Given two collections of $B-B$-bimodule maps $\{\phi_n : \otimes_B^n X
\rightarrow B\}_{n \in {\mathbb N}}$ and $\{\kappa_n : \otimes_B^n X
\rightarrow B\}_{n \in {\mathbb N}}$ extended multiplicatively, the
  following conditions are all equivalent:
\begin{itemize}
\item[(1)] $\phi_n = \sum_{\pi \in NC(n)} \kappa_\pi$ for each $n \in {\mathbb N}$.
\item[(2)] $\kappa_n = \sum_{\pi \in NC(n)} \mu(\pi,1_n) \phi_\pi$ for each $n \in {\mathbb N}$.
\item[(3)] $\phi_\tau = \sum_{\pi \in NC(n), \pi \leq \tau} \kappa_\pi$ for each $n \in {\mathbb N}, \tau \in NC(n)$.
\item[(4)] $\kappa_\tau = \sum_{\pi \in NC(n), \pi \leq \tau} \mu(\pi,\tau) \phi_\pi$ for each $n \in {\mathbb N}, \tau \in NC(n)$.
\end{itemize}
\end{proposition}

\begin{proof}[Sketch of Proof] Clearly $(3) \Rightarrow (1)$ and $(4) \Rightarrow (2)$ by taking $\tau = 1_n$. Next, suppose $(2)$ is given. 
We will prove $(4)$ by induction on the number of classes of $\tau$.
The basis case when $\tau = 1_n$ is clearly true.
If, on the other hand, $\tau = \rho \cup 1_{[k+1,l]}$ for $\rho \in NC([1,n] \setminus [k+1,l])$, we compute
\begin{eqnarray*}
\lefteqn{\kappa_\tau(x^1\otimes x^2\otimes \cdots\otimes x^n)}\\
 &=&\kappa_\rho(x^1\otimes \cdots\otimes x^{k-1}\otimes x^{k}\kappa_{l-k}(x^{k+1}\otimes \cdots\otimes x^l)\otimes x^{l+1}\otimes \cdots\otimes x^n)\\
               &=&  \sum_{\stackrel{{\lambda \in NC([1,n] \setminus [k+1,l]), \lambda \leq \rho}}{\nu \in NC([k+1,l])}} \mu(\lambda,\rho)\mu(\nu,1_{l-k}) \times \\ &&\phi_\lambda(x^1\otimes \cdots\otimes x^{k-1}\otimes x^{k}\phi_\nu(x^{k+1}\otimes \cdots\otimes x^l)\otimes x^{l+1}\otimes \cdots\otimes x^n)\\
               &=& \sum_{\pi \in NC(n), \pi \leq \tau} \mu(\pi,\tau)
               \phi_\pi(x^1\otimes  x^2\otimes \cdots\otimes x^n).
\end{eqnarray*}
Here, the first equality is by the multiplicativity of $\kappa$;
the second follows by (two applications of) the inductive
assumption; and the third equality follows from (i) the identification
$[0_n,\tau] = [0_{[1,n] \setminus [k+1,l]},\rho] \times [0_{[k+1,l]},1_{[k+1,l]}]$
of posets, (ii) the fact that $\mu$ is `multiplicative' with respect
to such decompositions of `intervals'; and from (iii) the
multiplicativity of $\phi$. This
finishes the inductive step and hence proves $(4)$.
An even easier proof shows that $(1) \Leftrightarrow (3)$.
Finally, $(3) \Leftrightarrow (4)$ by usual Mobius inversion in the
poset $NC(n)$.
\end{proof}

\begin{definition}
The free cumulants of a $B$-valued probability space $(A,\phi)$
are the $B-B$-bimodule maps $\kappa_n : \otimes_B^n A \rightarrow B$ 
associated as in Proposition \ref{Mobius} to the collection of $B-B$-bimodule maps
 $\{\phi_n : \otimes_B^n A \rightarrow B\}_{n \in {\mathbb N}}$ defined by $\phi_n(a^1\otimes \cdots \otimes a^n) = \phi(a^1a^2\cdots a^n)$. 
\end{definition}

The importance of the operator valued free cumulants lies in the following
theorem of Speicher linking their vanishing to freeness with amalgamation
over the base.

\begin{theorem}\label{cumulant}
Let $(A,\phi)$ be a $B$-valued probability space and $\{A_i : i \in I\}$ be a family of $B$-subspaces of $A$
such that $A_i$ is generated as an algebra over $B$ by $G_i \subseteq A_i$. This family is freely independent with amalgamation over $B$
iff for each positive integer $k$, indices $i_1,\cdots,i_k \in I$ that
are not all equal and elements $a_t \in G_{i_t}$ for $t = 1,2,\cdots,k$,
the equality $\kappa_k(a_1 \otimes a_2 \otimes \cdots \otimes a_k) = 0$ holds. \qed
\end{theorem}

We will regard $P_0(\Gamma,0) \subseteq Gr(\Gamma,0)$ as an
operator valued probability space with the map $\phi$ being defined on
$P_n(\Gamma,0)$ by the sum of the action of all $Hom([n],[0])$
morphisms. Equivalently, it
is the transport to the $Gr(\Gamma,0)$ picture
of the map given in the $F(\Gamma,0)$ picture by the `orthogonal projection
to $P_0(\Gamma,0)$'. 
This is easily checked to be an identity preserving $P_0(\Gamma,0)-P_0(\Gamma,0)$-bimodule map. Further, it preserves the faithful,
positive trace $\tau$, as is checked by definiton of $t$ in the $F(\Gamma,0)$ picture, and transporting to $Gr(\Gamma,0)$.

In order to state the main result of this section, we need to introduce
some notation. Observe first that $Gr(\Gamma,0)$ is generated
as an algebra by $P_0(\Gamma,0)$
and all $[\xi]$ where $\xi$ is a path length 2 in $\Gamma$.
For any odd vertex $w \in V_1$ (resp. even vertex $v \in V_0$), let
$\Gamma_w$ (resp. $\Gamma_v$) denote the subgraph of $\Gamma$ induced
on the vertex set  $V_0 \cup \{ w \}$ (resp., $\{v\} \cup V_1$).
The $\mu$ function of $\Gamma_w$
(resp. $\Gamma_v$) is the (appropriately normalised) restricted $\mu$ function of $\Gamma$.
Then, $Gr(\Gamma_w,0)$ is naturally isomorphic - as a $*$-probability space - to
the subalgebra of 
$Gr(\Gamma,0)$ generated by $P_0(\Gamma,0)$ and all $[\xi]$ such that $l(\xi) =2$ and $v^\xi_1 = w$, i.e., paths of length 2 with middle vertex $w$. 
We will refer to this subalgebra as $Gr(\Gamma_w,0)$.

Our main result in this section is then the following proposition.

\begin{proposition}\label{freeness} For an arbitrary finite, weighted, bipartite graph $\Gamma$, we have
\[Gr(\Gamma,0) = \ast_{P_0(\Gamma,0)} \{ Gr(\Gamma_w,0) : w \in V_1\} .\]
\end{proposition}

The crucial step in the proof of this proposition is the
identification of the $P_0(\Gamma,0)$-valued free cumulants
on the generators $[\xi]$ of
$Gr(\Gamma,0)$, which is done in the next proposition.  
In this, we will use a natural bijection 
$\pi \mapsto S(\pi)$ between the sets
$NC(\Sigma)$ and $TL(\Sigma \times \{1,2\})$ (for a totally ordered finite set $\Sigma$ and where we consider the
dictionary order on $\Sigma \times \{1,2\}$) defined as follows.
Suppose that $\pi \in NC(\Sigma)$. Let $C$ be a class of $\pi$
and enumerate the elements of $C$ in increasing order as, say, $C =
\{c_1,c_2,\cdots,c_t\}$. Decree $\{(c_1,2),(c_2,1)\}$, $\{(c_2,2),(c_3,1)\}$,
$\cdots$, $\{(c_{t-1},2),(c_t,1)\}$, $\{(c_t,2),(c_1,1)\}$ to be classes of $S(\pi)$.
Do this for each class of $\pi$ to define $S(\pi)$. 
Observe that $S(\pi)$ is a union of equivalence relations on $C \times \{0,1\}$ as $C$ varies over classes of $S$.
If $\Sigma = \{1,2,\cdots,n\}$, we will regard $S(\pi)$ as an element of  $TL(\{1,2,\cdots,2n\})$ or equivalently as an an element of $Hom([2n],[0])$, via the obvious
order isomorphisms.
We illustrate with an example. Suppose that $\pi = \{\{1,6\},\{2,3,4,5\}\}$.
\begin{figure}[htbp]
\psfrag{1}{\tiny $(1,1)$}
\psfrag{2}{\tiny $(1,2)$}
\psfrag{3}{\tiny $(2,1)$}
\psfrag{4}{\tiny $(2,2)$}
\psfrag{5}{\tiny $(3,1)$}
\psfrag{6}{\tiny $(3,2)$}
\psfrag{7}{\tiny $(4,1)$}
\psfrag{8}{\tiny $(4,2)$}
\psfrag{9}{\tiny $(5,1)$}
\psfrag{10}{\tiny $(5,2)$}
\psfrag{11}{\tiny $(6,1)$}
\psfrag{12}{\tiny $(6,2)$}
\begin{center}
\includegraphics[width=3in]{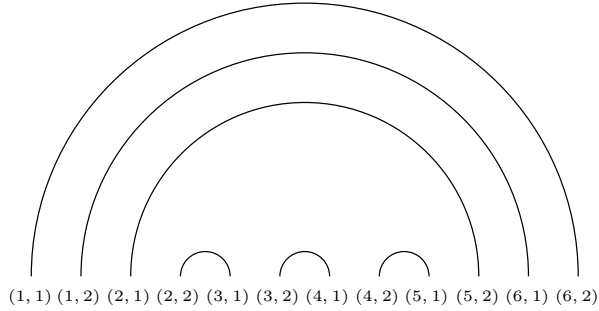}
\caption{The Temperley-Lieb equivalence relation $S(\pi)$}
\label{spifig}
\end{center}
\end{figure}
Then, $S(\pi)$ is shown in Figure \ref{spifig}. Regarded as
a Temperley-Lieb relation on $\{1,2,\cdots,12\}$, $S(\pi) = \{\{1,12\},\{2,11\},\{3,10\},\{4,5\},\{6,7\},\{8,9\}\}$.

We will also use the notion of 
a starry path in $\Gamma$ by which we
mean an even length path $\xi = (v^\xi_0 \stackrel{\xi_1}{\rightarrow}
v^\xi_1 \stackrel{\xi_2}{\rightarrow} v^\xi_2
\stackrel{\xi_3}{\rightarrow} \cdots \stackrel{\xi_{2n}}{\rightarrow}
v^\xi_{2n})$, where $\xi_{2i} = \widetilde{\xi_{2i+1}}$ (indices
modulo $2n$) for $i = 1,2, \cdots,n$. In such a path, we have $v^\xi_0 =
v^\xi_{2n}$ and all the odd $v^\xi_i$ are equal (to the centre of
the `star').
Note that a path $\xi$ of length $2n$ is starry exactly when
$S(1_n)([\xi]) \neq 0$ (in which case, it is a scalar multiple
of $e_v$ where $v$ is the start and finish point of $\xi$.)

\begin{proposition}\label{cumulants}
Let $\xi^1,\xi^2,\cdots,\xi^n$ be paths
in $\Gamma$ of length 2. The $P_0(\Gamma,0)$-valued
free cumulant on $Gr(\Gamma,0)$ is given thus:
$\kappa_n([\xi^1],[\xi^2],\cdots,[\xi^n])$ is non-zero only when $\xi = 
\xi^1 \circ \xi^2 \circ \cdots \circ \xi^n$ is a starry path in
$\Gamma$  (in particular this composition should make sense),
in which case, all the odd $v^\xi_i$ are equal to some $w \in V_1$ and
$v^\xi_0 = v^\xi_{2n}$ is some $v \in V_0$, and 
\begin{eqnarray*}
\kappa_n([\xi^1],[\xi^2],\cdots,[\xi^n]) &=& S(1_n)([\xi]) \\
&=& \frac
{
\mu(v^\xi_2)\mu(v^\xi_4)\cdots \mu(v^\xi_{2n-2})
}
{
\mu(w)^{n-2}\mu(v)
} e_v.
\end{eqnarray*}
\end{proposition}

\begin{proof}
Consider the $P_0(\Gamma,0)-P_0(\Gamma,0)$-bimodule $P_2(\Gamma,0)$ and define for each $n \geq 1$, $\widetilde{\kappa}_n : P_2(\Gamma,0) \times P_2(\Gamma,0)
 \times \cdots \times P_2(\Gamma,0) \rightarrow P_0(\Gamma,0)$ by the
$\C$-multilinear extension of the prescription given on basis
elements by the statement of the proposition. It is easily checked that
$\tilde{\kappa}_n$ induces a bimodule map also denoted
$\tilde{\kappa}_n : \otimes_{P_0(\Gamma,0)}^n P_2(\Gamma,0) \rightarrow
P_0(\Gamma,0)$.

We will check that $\widetilde{\kappa}_n$ agrees on
$\otimes_{P_0(\Gamma,0)}^n P_2(\Gamma,0)$ with the operator valued
cumulant $\kappa_n$.
In view of Proposition \ref{Mobius}, it suffices to check that if $\widetilde{\kappa}_\pi$ is the multiplicative extension of $\widetilde{\kappa}$, then,
$$
\phi([\xi^1 ] \bullet  [\xi^2 ] \bullet  \cdots  \bullet  [\xi^n]) =
\sum_{\pi \in NC(n)} \widetilde{\kappa}_\pi([\xi^1] \otimes \cdots
\otimes [\xi^n]) 
$$
for all paths $\xi^1,\cdots,\xi^n$ of length 2 in $\Gamma$.

Observe first that both sides of the equation above vanish unless
$\xi^1 \circ \xi^2 \circ \cdots \circ \xi^n$ makes sense and defines
a path $\xi$ with equal end-points, so we suppose this to be the case. 
By definition then, the left hand side $\phi([\xi^1 ] \bullet  [\xi^2 ] \bullet  \cdots  \bullet  [\xi^n])$
of the desired equality is given by $\sum_{S \in Hom([2n],[0])} S([\xi^1 \circ \xi^2 \circ \cdots \circ \xi^n])$.
In view of the natural bijection between $NC(n)$ and $Hom([2n],[0])$
alluded to above, it clearly suffices to see that
\begin{equation}\label{stp} \widetilde{\kappa}_\pi([\xi^1] \otimes
  \cdots \otimes [\xi^n]) = S(\pi)([\xi^1 \circ  \xi^2 \circ \cdots
  \circ \xi^n]).\end{equation} 
  
Observe that for $\pi \in NC(n)$ and paths $\xi^1,\xi^2,\cdots,\xi^n$
of length $2$ for which the composite $\xi^1 \circ \cdots \circ \xi^n$
is defined, 
\begin{eqnarray}\label{spi}
\lefteqn{S(\pi)([\xi^1 \circ  \xi^2 \circ \cdots
\circ \xi^n]) =} \nonumber \\
&\left(  \prod_{C = \{c_1,\cdots,c_t\} \in \pi}\left\{  \delta_{\xi^{c_t}_2,\widetilde{\xi^{c_{1}}_1}} \frac{\mu(v^{\xi^{c_{1}}}_1)}{\mu(v^{\xi^{c_t}}_2)} \prod_{p=1}^{t-1} \left( \delta_{\xi^{c_p}_2,\widetilde{\xi^{c_{p+1}}_1}} \frac{\mu(v^{\xi^{c_p}}_2)}{\mu(v^{\xi^{c_{p+1}}}_1)} \right)\right\}  \right) e_{v^{\xi^1}_0}
\end{eqnarray}
(where $c_1 < c_2 < \cdots < c_t$).
For instance, for the $S$ shown in Figure \ref{spifig}, we have
\begin{eqnarray*}
S(\pi)([\xi^1 \circ  \xi^2 \circ \cdots
\circ \xi^n]) &=&
\left\{
\delta_{\xi^6_2,\widetilde{\xi^1_1}} \frac{\mu(v^{\xi^{{1}}}_1)}{\mu(v^{\xi^{6}}_2)}
\delta_{\xi^1_2,\widetilde{\xi^6_1}} \frac{\mu(v^{\xi^{{1}}}_2)}{\mu(v^{\xi^{6}}_1)} \right\} \times\\
&&
\left\{
\delta_{\xi^5_2,\widetilde{\xi^2_1}} \frac{\mu(v^{\xi^{{2}}}_1)}{\mu(v^{\xi^{5}}_2)}
\delta_{\xi^2_2,\widetilde{\xi^3_1}} \frac{\mu(v^{\xi^{{2}}}_2)}{\mu(v^{\xi^{3}}_1)}
\delta_{\xi^3_2,\widetilde{\xi^4_1}} \frac{\mu(v^{\xi^{{3}}}_2)}{\mu(v^{\xi^{4}}_1)}
\delta_{\xi^4_2,\widetilde{\xi^5_1}} \frac{\mu(v^{\xi^{{4}}}_2)}{\mu(v^{\xi^{5}}_1)}
 \right\}
 e_{v^{\xi^1}_0}.
\end{eqnarray*}

We prove equation (\ref{stp})  by induction on the number of classes
of $\pi$. In case $\pi = 1_n$, 
it holds by definitions of $\widetilde{\kappa}_\pi$ and $S$.
Suppose next that $\pi = \rho \cup 1_{[k+1,l]}$ for $\rho \in NC([1,n]
\setminus [k+1,l])$. Since $\widetilde{\kappa}$ multiplicatively
extends $\{\widetilde{\kappa_n}\}_{n \in {\mathbb N}}$ we have, 
$$
\widetilde{\kappa}_\pi([\xi^1],\cdots,[\xi^n]) =
\widetilde{\kappa}_\rho([\xi^1],\cdots,[\xi^{k-1}],[\xi^k] \widetilde{\kappa}_{l-k}([\xi^{k+1}],\cdots,[\xi^l]),[\xi^{l+1}],\cdots,[\xi^n]).
$$
By definition,
$$
\widetilde{\kappa}_{l-k}([\xi^{k+1}],\cdots,[\xi^l]) =\left\{ \delta_{\xi^{l}_2,\widetilde{\xi^{k+1}_1}} \frac{\mu(v^{\xi^{k+1}}_1)}{\mu(v^{\xi^{l}}_2)} \prod_{p=1}^{l-k-1} \left( \delta_{\xi^{p+k}_2,\widetilde{\xi^{{p+k+1}}_1}} \frac{\mu(v^{\xi^{p+k}}_2)}{\mu(v^{\xi^{{p+k+1}}}_1)} \right)\right\} e_{v^{\xi^{k+1}}_0}.
$$
Since $f(\xi^k) = s(\xi^{k+1}) = v^{\xi^{k+1}}_0$, we conclude that
\begin{eqnarray*}
\lefteqn{\widetilde{\kappa}_\pi([\xi^1],\cdots,[\xi^n])}\\
&=&
\left\{ \delta_{\xi^{l}_2,\widetilde{\xi^{k+1}_1}}
  \frac{\mu(v^{\xi^{k+1}}_1)}{\mu(v^{\xi^{l}}_2)} \prod_{p=1}^{l-k-1}
  \left( \delta_{\xi^{p+k}_2,\widetilde{\xi^{{p+k+1}}_1}}
    \frac{\mu(v^{\xi^{p+k}}_2)}{\mu(v^{\xi^{{p+k+1}}}_1)} \right)
\right\}\\
&& ~~~~~~~~~~~~~~~~~~~~~~~~~~\times ~~~~\widetilde{\kappa}_\rho([\xi^1],\cdots,[\xi^{k-1}],[\xi^k] ,[\xi^{l+1}],\cdots,[\xi^n])\\
&=& \left\{ \delta_{\xi^{l}_2,\widetilde{\xi^{k+1}_1}}
  \frac{\mu(v^{\xi^{k+1}}_1)}{\mu(v^{\xi^{l}}_2)} \prod_{p=1}^{l-k-1}
  \left( \delta_{\xi^{p+k}_2,\widetilde{\xi^{{p+k+1}}_1}}
    \frac{\mu(v^{\xi^{p+k}}_2)}{\mu(v^{\xi^{{p+k+1}}}_1)} \right)
\right\} \\
&& ~~~~~~~~~~~~~~~~~ \times\left( \prod_{C = \{c_1,\cdots,c_t\} \in \rho}\left\{  \delta_{\xi^{c_t}_2,\widetilde{\xi^{c_{1}}_1}} \frac{\mu(v^{\xi^{c_{1}}}_1)}{\mu(v^{\xi^{c_t}}_2)} \prod_{p=1}^{t-1} \left( \delta_{\xi^{c_p}_2,\widetilde{\xi^{c_{p+1}}_1}} \frac{\mu(v^{\xi^{c_p}}_2)}{\mu(v^{\xi^{c_{p+1}}}_1)} \right)\right\} \right)  e_{v^{\xi^1}_0}\\
&=&\left(  \prod_{C = \{c_1,\cdots,c_t\} \in \pi}\left\{  \delta_{\xi^{c_t}_2,\widetilde{\xi^{c_{1}}_1}} \frac{\mu(v^{\xi^{c_{1}}}_1)}{\mu(v^{\xi^{c_t}}_2)} \prod_{p=1}^{t-1} \left( \delta_{\xi^{c_p}_2,\widetilde{\xi^{c_{p+1}}_1}} \frac{\mu(v^{\xi^{c_p}}_2)}{\mu(v^{\xi^{c_{p+1}}}_1)} \right)\right\}  
\right) e_{v^{\xi^1}_0}\\
&=& S(\pi)([\xi^1 \circ  \xi^2 \circ \cdots
\circ \xi^n])
\end{eqnarray*}
where the second equality follows from the inductive assumption
and equation (\ref{spi}) applied to $\rho$,
and the last equality follows from equation (\ref{spi})
applied to $\pi$.
\end{proof}

The proof of Proposition \ref{freeness} is now immediate.

\begin{proof}[Proof of Proposition \ref{freeness}] Since $Gr(\Gamma,0)$ is clearly generated by all the $Gr(\Gamma_w,0), w \in V_1$, it needs only
to be seen that the family $\{Gr(\Gamma_w,0) : w \in V_1\}$ is free with amalgamation over $P_0(\Gamma,0)$. This follows from Theorem \ref{cumulant} and Proposition
\ref{cumulants}.
\end{proof}

Let $\lambda : Gr(\Gamma,0) \rightarrow {\mathcal L}(H(\Gamma,0))$
also denote the composite of the isomorphism of $Gr(\Gamma)$ with
$F(\Gamma)$ and $\lambda : F(\Gamma,0) \rightarrow {\mathcal L}(H(\Gamma,0))$.
Thus, $M(\Gamma,0) = \lambda(Gr(\Gamma,0))^{\prime\prime}$.
It now follows fairly easily - see Proposition 4.6 of \cite{KdySnd2008},
for instance - that $M(\Gamma_w,0) \cong \lambda(Gr(\Gamma_w,0))^{\prime\prime}$, the content in this statement being that $Gr(\Gamma_w,0)$ is interpreted as a subalgebra of $Gr(\Gamma,0)$.
We will thus identify $M(\Gamma_w,0)$ with $\lambda(Gr(\Gamma_w,0))^{\prime\prime} \subseteq M(\Gamma,0)$.

Now, by general principles, Proposition \ref{freeness} extends to its
von Neumann completions - meaning that
\begin{equation}\label{mgamam}
M(\Gamma,0) = \ast_{P_0(\Gamma,0)} \{ M(\Gamma_w,0) : w \in V_1\},
\end{equation}
and similarly, by interchanging the roles of 0 and 1, we have
\begin{equation}\label{mgamam1}
M(\Gamma,1) = \ast_{P_0(\Gamma,1)} \{ M(\Gamma_v,1) : v \in V_0\}.
\end{equation}

\section{Graphs with a single odd vertex}

In this section we fix the following notation. Let $\Lambda$ be a
graph with at least one edge and a single odd vertex $w$ and even vertices $v_1,\cdots,v_l$. We assume that for $i = 1,\cdots, k$, the vertex $v_i$
is joined to $w$ by $q_i > 0$ edges, while the vertices $v_i$
for $i = k+1,\cdots,l$ are isolated.  Thus $k \geq 1$.
We also set $\mu^2(v_i) = \alpha_i$ and $\mu^2(w) = \beta$, so that
$\beta + \sum_{i=1}^l \alpha_i = 1$. 
Our goal in this section is the
explicit determination of the finite von Neumann algebra $M(\Lambda,0)$.

We begin with a simple observation. The assignments of
$Gr(\Gamma)$ or $M(\Gamma)$ to a graph $\Gamma$ clearly take
disjoint unions to (appropriately weighted) direct sums. Thus, if
$\tilde{\Lambda}$ denotes the connected component of $w$ in the graph $\Lambda$ (with normalised restricted $\mu$), then
\begin{equation}\label{ll'}
M(\Lambda) \cong \underset{\gamma}{M(\tilde{\Lambda})} \oplus
\underset{\alpha_{k+1}}{\mathbb C} \oplus
\underset{\alpha_{k+2}}{\mathbb C} \oplus \cdots
\oplus
\underset{\alpha_l}{\mathbb C}
\end{equation}
where $\gamma = 1-\sum_{i=k+1}^l \alpha_i$.
Note that $M(\Lambda,1) = M(\tilde{\Lambda},1)$. We begin by analysing $M(\Lambda,1)$ when
$k=1$.

\begin{remark}\label{lf}
We shall adopt the convention that $LF(1) = L{\mathbb Z}$ whereas
$LF(r)$, for $1 < r < \infty$, will be referred to as an interpolated free group factor with finite parameter.
\end{remark}

\begin{proposition}\label{line}
Let $\Omega$ be a graph with a single even vertex $v$ and single odd vertex $w$ joined by $q>0$ edges. Let $\mu^2(v) = \alpha$ and
$\mu^2(w) = \beta = 1-\alpha$. Then,
$$
M(\Omega,1) \cong 
\left\{ 
\begin{array}{ll}
                   LF(q^2) & {\text {if\ }} \frac{\alpha}{\beta} > q \\
                    LF(\frac{2q\alpha}{\beta} - \frac{\alpha^2}{\beta^2})  & {\text {if\ }} \frac{1}{q} \leq \frac{\alpha}{\beta} \leq q\\
                    \underset{1-\frac{\alpha q}{\beta}}{\mathbb C} \oplus
\underset{\frac{\alpha q}{\beta}}{LF(2-\frac{1}{q^2}})  & {\text {if\ }} \frac{\alpha}{\beta} < \frac{1}{q}
        \end{array}
\right. 
$$
\end{proposition}

\begin{proof}
Consider $Gr(\Omega,1)$ generated by $q^2$ paths of length 2 based at $w$. Denoting the path $(w \stackrel{i}{\rightarrow} v \stackrel{j}{\rightarrow} w)$ by $e_{ij}$, the operator-valued (in this
 case, scalar valued) free cumulant calculation of Proposition \ref{cumulants} implies that  the
$q \times q$-matrix $X = (([e_{ij}]))$ is a uniformly $R$-cyclic matrix
- in the sense of Definition 10 of \cite{KdySnd2009} - with determining
sequence $\alpha_t = \left(\frac{\mu(w)}{\mu(v)}\right)^{t-2}$.
Theorem 11 of \cite{KdySnd2009} now implies that $X$ is free
Poisson with rate $\frac{\alpha}{\beta q}$.
Now, the proof of Proposition 24 of \cite{KdySnd2009} may be imitated to yield the desired result.
\end{proof}

As a consequence, we single out a crisp determination of precisely when $M(\Omega)$ is a factor.

\begin{corollary}
Let $\Omega$ be as in Proposition \ref{line}. Then,
$M(\Omega)$ is a factor if and only if $q > 1$ and $\frac{1}{q} \leq \frac{\alpha}{\beta} \leq q$, in which case $M(\Omega) \cong LF(1+2q\alpha\beta-\alpha^2-\beta^2)$.
\end{corollary}

\begin{proof} If $M(\Omega)$ is a factor, then, so is $M(\Omega,1)$ and it follows from Proposition \ref{line}
that $\frac{1}{q} \leq \frac{\alpha}{\beta}$.
Similarly the factoriality of $M(\Omega,0)$ and Proposition \ref{line}
will imply that $\frac{1}{q} \leq \frac{\beta}{\alpha}$.
Thus, 
 $\frac{1}{q} \leq \frac{\alpha}{\beta} \leq q$. 
 To see that $q > 1$, it suffices to observe that if $q=1$, the already proved inequality shows that $\alpha = \beta$ and 
 then again by Proposition \ref{line}, both $M(\Omega,0)$
and $M(\Omega,1)$ are $LF(1) \cong L{\mathbb Z}$ so $M(\Omega)$
cannot be a factor.

For the converse, if $q > 1$ and $\frac{1}{q} \leq \frac{\alpha}{\beta} \leq q$, then $q \neq \frac{1}{q}$ and so at least one of the inequalities
among $\frac{\alpha}{\beta} \geq \frac{1}{q}$ and $q \geq  \frac{\alpha}{\beta}$
must be strict. Hence 
\begin{eqnarray*}
\frac{2q\alpha}{\beta} - \frac{\alpha^2}{\beta^2} &=& \frac{\alpha}{\beta}(2q-\frac{\alpha}{\beta})\\
&>& \frac{1}{q}.q\\
&=& 1,
\end{eqnarray*}
and so $M(\Gamma,1)$ is an  interpolated free group factor with finite parameter. 
Similarly, so is $M(\Gamma,0)$. By Lemma \ref{locglob}, $M(\Gamma)$ is a factor. Now the corner
formula for interpolated free group factors - see \cite{Dyk1994} or \cite{Rdl1994} - implies that
$M(\Omega) \cong LF(1+2q\alpha\beta-\alpha^2-\beta^2)$.
\end{proof}

We now wish to analyse $M(\Lambda,1) = M(\tilde{\Lambda},1)$. For this, recall
that the weighting, say $\tilde{\mu}$, on $\tilde{\Lambda}$ is given
by the normalised restriction of $\mu$. Thus $\tilde{\mu}^2(v_i) =a_i$ for $1 \leq i \leq k$ and $\tilde{\mu}^2(w) =b$, where $\sum_{i=1}^k a_i + b = 1$ and $(\alpha_1:\cdots:\alpha_k:\beta) = (a_1:\cdots:a_k:b)$.

\begin{proposition}\label{ml1} With the foregoing notation,
$$
M(\tilde{\Lambda},1) \cong
\left\{
\begin{array}{ll}
LF\left(\sum_{\{i: q_ib < a_i\}} q_i^2 + \sum_{\{i : q_ib \geq a_i\}} \frac{2q_ia_i}{b} - \left(\frac{a_i}{b}\right)^2\right) & {\text {if\ }} b \leq \sum_{i=1}^k q_ia_i\\
\underset{1-\sum_i \frac{q_ia_i}{b}}{\mathbb C} \oplus
\underset{\sum_i \frac{q_ia_i}{b}}{LF(2-\frac{\sum_i a_i^2}{(\sum_i q_ia_i)^2}})  & {\text {if\ }} b > \sum_{i=1}^k q_ia_i
\end{array}
\right.
$$
In particular, $M(\tilde{\Lambda},1)$ is an $LF(r)$ for some $r \geq 1$ iff $b \leq \sum_{i=1}^k q_ia_i$. 
\end{proposition}

\begin{proof} 
By equation (\ref{mgamam1}), we have
$$
M(\tilde{\Lambda},1) \cong *_{P_0(\tilde{\Lambda},1)} \{ M(\tilde{\Lambda}_{v_i},1) : i =1,2,\cdots,k \} \cong *_{\ }  \{ M(\tilde{\Lambda}_{v_i},1) : i =1,2,\cdots,k \},
$$
where the second isomorphism holds since $P_0(\tilde{\Lambda},1) \cong {\mathbb C}$. 
Each $M(\tilde{\Lambda}_{v_i},1)$ is determined using Proposition \ref{line}. Now computations from [Dyk1993] - see Proposition 1.7 - and a little calculation finish the proof.
\end{proof}

\begin{proposition}\label{lf+fd}
If $\Lambda$  has a single odd vertex and at least two edges, then $M(\Lambda) \cong LF(s) \oplus A$, for some finite $s > 1$ and a finite-dimensional abelian  $A$.
\end{proposition}

\begin{proof}
Notice that $\tilde{\Lambda}$ satisfies the hypotheses of this proposition and in addition, is connected. 
In view of equation (\ref{ll'}), it suffices to prove the proposition for
$\tilde{\Lambda}$; in other words, we may assume without loss of
generality that $\Lambda$ is connected.

Hence Theorem \ref{Mgamstr} is applicable and  $M(\Lambda)$ has the form $M \oplus A$ for some $II_1$ factor $M$ and a finite-dimensional abelian $A$.

Now Proposition \ref{ml1} tells us that some corner of $M(\Lambda,1)$ and hence of $M(\Lambda)$ is an $LF(r)$  for some finite $r$. 
On the other hand,   the hypothesis that $\Lambda$ has at least two edges ensures that $M(\Lambda,1)$ is not commutative and hence $r > 1$.
This corner is necessarily a corner of $M \cong LF(s)$ for some finite $s >1$.
\end{proof}

\begin{corollary}\label{l0str}
Let $\Lambda$ be any graph with a single odd vertex and non-empty edge set $E$. Then,
\begin{equation*}
M(\Lambda,0) \cong \left\{ \begin{array}{ll}
                                           LF(s) \oplus A &\mbox{ if } |E| > 1\\
                                           L{\mathbb Z} \oplus A &\mbox{ if } |E| = 1.
                                           \end{array}
                                 \right.
\end{equation*}        
for some $s>1$ and finite-dimensional abelian $A$.  
\end{corollary}

\begin{proof}
In case $|E|>1$, $M(\Lambda,0)$ is necessarily non-abelian and the
desired assertion  is a direct consequence of Proposition \ref{lf+fd}.
When $|E| =1$, observe, as in equation (\ref{ll'}), that $M(\Lambda,0) = M(\tilde{\Lambda},0) \oplus A$  for some finite-dimensional abelian $A$. Now the desired result follows from Proposition \ref{line} applied
with $\Omega$ being $\tilde{\Lambda}$ with vertex parity reversed.
(This is because the $q$ of Proposition \ref{line} is $1$ and the parameter occurring in the $LF(\cdot)$ factor is 1 in all the three cases considered there.)
\end{proof}

\section{The structure of $M(\Gamma)$}

In this section, we determine the structure of $M(\Gamma)$
for any finite, connected, bipartite graph $\Gamma$ with Perron-Frobenius weighting. The main technical result used in the proof
is Theorem \ref{amal} which is a consequence of the results in
\cite{Dyk2009}.

\begin{theorem}\label{amal}
Let $M(w), w \in V_1$ be a finite family of tracial von Neumann algebraic
probability spaces over a finite-dimensional abelian probability space
$D$.
Suppose that each $M(w) \cong LF(r_w) \oplus A(w)$ 
with $1 \leq r_w < \infty$ and finite-dimensional abelian $A(w)$
 and that $M = *_D \{M(w) : w \in V_1\}$ is a factor.
Then, $M$ is either an
interpolated
free group factor with finite parameter or the hyperfinite ($II_1$) factor.
\end{theorem}

The following theorem is one of the main results of this paper.

\begin {theorem}\label{mgamlf}
Let a connected graph $\Gamma$  with at least two edges be
equipped with the Perron-Frobenius weighting. Then $M(\Gamma) \cong LF(s)$ for some $1 < s < \infty$.
\end{theorem}

\begin{proof}
By Corollary \ref{mgamfac}, $M(\Gamma)$ is a $II_1$-factor and so, to see
that it an interpolated free group factor with finite parameter, it suffices to see that the corner $M(\Gamma,0)$ is also one.

The hypotheses on $\Gamma$ ensure that 
Corollary \ref{l0str} is applicable to
$\Gamma_w$ 
for each odd vertex $w$.
Then it follows from equation (\ref{mgamam}) that  the hypotheses of Theorem \ref{amal} are satisfied 
with $D = P_0(\Gamma,0)$, $M(w) = M(\Gamma_w,0)$ 
for $w \in V_1$ and $M = M(\Gamma,0)$ and so $M(\Gamma,0)$ is
either an interpolated free group factor with finite parameter or the hyperfinite factor.

To conclude the proof, we only need to ensure that $M(\Gamma,0)$
is not hyperfinite. For this we consider two cases.
\begin{enumerate}
\item[Case 1.] Suppose some odd vertex $w$ of $\Gamma$ has degree greater
than 1.
In this case Corollary \ref{l0str} shows that $LF(r_w)$ for some $r_w >1$ 
is a corner of $M(\Gamma_w,0)$. A corner of a subalgebra of 
the hyperfinite factor cannot be $LF(r_w)$ (which is not injective).
Hence $M(\Gamma,0)$ is not hyperfinite.
\item[Case 2.] Every odd vertex of $\Gamma$ has degree 1. Thus
$\Gamma$ is the complete bipartite graph $K(1,n)$ for $n \geq 2$.
The Perron-Frobenius weighting on this graph assigns $\frac{1}{\sqrt{2}}$
to the odd vertex and $\frac{1}{\sqrt{2n}}$ to each even vertex.
Now, for any odd vertex $w$, Proposition \ref{line} applied with $\Omega$ being $\Gamma_w$ with reversed vertex parity implies that
$M(\Gamma_w,0) \cong \underset{\ \ 1-\delta^{-1}}{{\mathbb C}} \oplus \underset{\ \ \delta^{-1}}{L{\mathbb Z}}$, where $\delta = \sqrt{n}$. Clearly, $P_0(\Gamma,0) \cong {\mathbb C}$. Therefore, from equation (\ref{mgamam}),
$$
M(\Gamma,0) \cong (\underset{\ \ 1-\delta^{-1}}{{\mathbb C}} \oplus \underset{\ \ \delta^{-1}}{L{\mathbb Z}})^{*n} \cong LF(2\sqrt{n}-1),
$$
where the second isomorphism is proved in Corollary 16 of \cite{KdySnd2009}.
Since $n\geq 2$, $M(\Gamma,0)$ is an interpolated free group factor
with finite parameter
in this case too.
\end{enumerate}
\end{proof}

The only connected graph to which Theorem \ref{mgamlf} does not apply is the graph $A_2$ which has two vertices joined by a single
edge. For completeness, we determine, in the following proposition,
the structure of $M(\Gamma)$ in this case.

\begin{proposition}
Let $\Gamma$ be the $A_2$ graph with a single even vertex $v$
and a single odd vertex $w$ joined by a single edge. Equip $\Gamma$
with its Perron-Frobenius weighting given by $\mu^2(v) = \frac{1}{2}
=\mu^2(w)$. Then $M(\Gamma) \cong M_2(L{\mathbb Z})$.
\end{proposition}

\begin{proof} Recall from \S 1 that with $\Gamma$ being the $A_2$
  graph, elements of $Gr(\Gamma)$ may be regarded as matrices with
  rows and columns indexed by the set $\{v,w\}$ and $(p,q)$ entry
  (with $p,q \in \{v,w\}$) being a linear combination of paths from
  $p$ to $q$.   We shall write $M_{ij}$ for $e_iM(\Gamma) e_j$ for $0 \leq i,j \leq
  1$, where of course $e_0=e_v$ (resp. $e_1=e_w$) denotes the projection onto
  the subspace $\CH_0$ (resp. $\CH_1$) of $\CH$ generated by the set
  of all paths starting at $v$ (resp. $w$). Let $\xi_n$
  (resp. $\eta_n$) be the unique path of length $n$ which starts at
  $v$ (resp. $w$). Then, $\CH = \CH_0 \oplus \CH_1$, and also (see
  equation (\ref{bracedef}))
  $\{\{\xi_n\}: n \geq 0\}$ (resp. $\{\{\eta_n\}: n \geq 0\}$) is an
  orthonormal basis for $\CH_0$ (resp. $\CH_1$).

Let $x = \lambda(\xi_1) \in M_{01}$. The definition of
multiplication in $F(\Gamma)$ shows that $x \{\xi_n\} = 0$ and $x \{\eta_n\} =
\{\eta_{n+1}\} + \{\eta_{n-1}\}$ for all $ n \geq 1$ (with $\{\eta_{-1}\}
= 0$). So, if $w:\CH_0 \rar \CH_1$ is the (obviously unitary) operator 
defined by $w(\{\xi_n\}) = \{\eta_n\}$, we see that $x=ws$ where $s
\in \CL(\CH_0)$ is the (standard semi-circular) operator given by
$s\{\xi_n\} = \{\xi_{n+1}\} + \{\xi_{n-1}\}$. It follows that $x$ is
injective (since the Wigner distribution has no atoms).

It follows that if $x = u|x|$ denotes the polar decomposition of $x$, then 
$u^*u = e_1$, and similarly one sees that $uu^*= e_0$.

Now if $y \in M_{01}$ is arbitrary, then $y=e_0ye_1$ and we see that
$yu^* = (e_0ye_1)(e_0ue_1)^*$ $= e_0ye_1u^*e_0 \in M_{00}$, and hence
$y= ye_1 = yu^*u \in M_{00}u$. Arguing similarly, we see that the
maps
\[a \mapsto au,  a \mapsto u^*a, \mbox{ and } a \mapsto u^*au\]
define linear isomorphisms of
$M_{00}$ onto $M_{01}, M_{10} \mbox{ and } M_{11}$
respectively. Finally, it is easy to see that the assignment
\[\left[ \ba{ll} a_{00} & a_{01}\\a_{10} & a_{11} \ea
  \right] \mapsto \left[ \ba{cc} a_{00} & a_{01}u\\u^*a_{10} & u^*a_{11}u \ea
  \right] \]
defines an isomorphism of $M_2(M_{00})$ onto $M(\Gamma)$.
Since $M_{00} \cong L\Z$ by Proposition~\ref{line}, the proof
is complete.
\end{proof}

\section{Application to the GJS construction}

In this section we relate our $Gr(\Gamma)$ to the sequence of algebras
$Gr_k(P)$ of \cite{GnnJnsShl2008}. Let $P$ be a subfactor planar
algebra with finite principal graph $\Gamma = (V,E)$, distinguished vertex $*$
and modulus $\delta > 1$. Thus $\delta$ is the Perron-Frobenius 
eigenvalue of $\Gamma$ and we let $\mu^2(\cdot)$ be the Perron-Frobenius
eigenvector normalised so that $\sum_{v \in V} \mu^2(v) = 1$.
Let $tr$ be the normalised picture trace on $P_n$.

Most of the following facts about the tower of algebras $$(P_{0_+} =) P_0 \subseteq P_1 \subseteq P_2 \subseteq \cdots $$ can all be found in
\cite{JnsSnd1997}.
For vertices $v,w$ of $\Gamma$, we write ${\mathcal P}_n(v,w)$ for the set of paths of length $n$ in $\Gamma$ beginning at $v$ and ending at $w$.
Similarly we use notation such as ${\mathcal P}_n(v,\cdot)$ for the set of paths of length $n$ in $\Gamma$ beginning at $v$.

$P_n$ has a basis given by pairs of paths $(\xi(+),\xi(-))$ in $\Gamma$
such that 
$\xi(\pm) \in {\mathcal P}_n(*,v)$ for some vertex $v \in V$.
The set ${\mathcal P}_{min}(Z(P_n))$ of minimal central projections of $P_n$ is in natural bijection with $\{v \in V : v = f(\xi) {\text {~for some~}} \xi \in {\mathcal P}_n(*,\cdot) \}$. For such a $v$,
denote the corresponding minimal central projection in $P_n$ by
$e(v,n)$ and any minimal projection under $e(v,n)$ by $p(v,n)$.
Then, $\{(\xi(+),\xi(-)) : \xi(\pm) \in {\mathcal P}_n(*,v) \}$ are matrix units 
(meaning $(\xi(+),\xi(-))(\eta(+),\eta(-)) = \delta^{\xi(-)}_{\eta(+)} (\xi(+),\eta(-))$ and $(\xi(+),\xi(-))^* = (\xi(-),\xi(+))$)
for the matrix algebra
$e(v,n)P_n$. Further, with $tr(\cdot)$ denoting the normalised picture trace on the planar algebra $P$, we have $tr(p(v,n)) = \delta^{-n} \frac{\mu^2(v)}{\mu^2(*)}$.

The inclusion of $P_n$ into $P_{n+1}$ is given by\\
\begin{eqnarray}\label{inclusion}
(\xi(+),\xi(-)) &\mapsto& 
\sum_{v \in V} \sum_{\rho(\pm) \in {\mathcal P}_{n+1}(*,v)}
\delta^{\xi(+)}_{\rho(+)_{[0,n]}}
\delta^{\xi(-)}_{\rho(-)_{[0,n]}}
\delta^{\rho(+)_{n+1}}_{\rho(-)_{n+1}}
(\rho(+),\rho(-))\\
&=& \sum_{\lambda \in {\mathcal P}_1(f(\xi(\pm),\cdot)}(\xi(+) \circ \lambda,\xi(-) \circ \lambda)\nonumber
\end{eqnarray}
The $\tau$-preserving conditional expectation $P_{n+1} \rightarrow P_n$ is given by
\begin{equation}\label{condexp}
(\xi(+),\xi(-)) \mapsto \delta^{\xi(+)_{n+1}}_{\xi(-)_{n+1}} 
\frac{\mu^2(v^{\xi(\pm)}_{n+1})}{\delta \mu^2(v^{\xi(\pm)}_{n})}
(\xi(+)_{[0,n]},\xi(-)_{[0,n]}).
\end{equation}
The Jones projection $e_n \in P_n$ for $n \geq 2$ is given by
\begin{equation}\label{jones}
\sum_{v \in V} \sum_{\xi(\pm) \in {\mathcal P}_n(*,v)}
 \delta^{\xi(+)_{[0,n-2]}}_{\xi(-)_{[0,n-2]}}
 \delta^{\xi(+)_{n-1}}_{\widetilde{\xi(+)_n}} \delta^{\xi(-)_{n-1}}_{\widetilde{\xi(-)_{n}}} \frac{\mu(v^{\xi(+)}_{n-1})\mu(v^{\xi(-)}_{n-1})}{\delta \mu^2(v^{\xi(\pm)}_n)}
(\xi(+),\xi(-)).
\end{equation}
In these formulae we have written $\delta^i_j$ for the Kronecker delta.

Our main observation is that the construction in \cite{GnnJnsShl2008} of $Gr_k(P)$
(after conjugating by suitable powers of the rotation tangle) depends
only on the actions of the inclusion, multiplication and right conditional expectation tangles - and not on that of the the left conditional
expectation tangles.
Hence, in principle, 
`$Gr_k(P)$ depends only on the graphs and not on the connection'.

\begin{figure}[htbp]
\psfrag{m}{ $\mapsto$}
\psfrag{2n}{\tiny $2n$}
\begin{center}
\includegraphics[height=4cm]{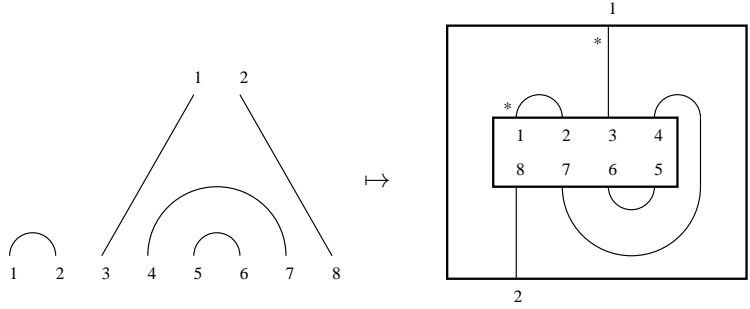}
\caption{Correspondence between ${\mathcal E}_{ev}$ and good Temperley-Lieb tangles}
\label{epitl}
\end{center}
\end{figure}

We first need to note that the action of the category epi-TL
or ${\mathcal E}$ of \S 1 on $Gr(\Gamma,*)$ is `essentially the same' as that of certain annular tangles on $Gr_0(P)$.
Consider the full subcategory of ${\mathcal E}$ consisting
only of the objects $[0],[2],[4],\cdots$. We will denote this
category by ${\mathcal E}_{ev}$. Any morphism in ${\mathcal E}_{ev}$, say an element of $Hom([2n],[2m])$, naturally yields an annular tangle
with an internal $n$-box and an external $m$-box as in the example
in Figure \ref{epitl} for $n=4,m=1$.
This identification of $Hom([2n],[2m])$, composed with the action of annular tangles on planar algebras is seen to yield an ${\mathcal E}_{ev}$ action on $\{P_n\}_{n \geq 0}$.
(The tangles in the image of ${\mathcal E}_{ev}$ are what were called
$0$-good annular tangles in \cite{KdySnd2008}.)

We will find it convenient to identify a basis element $(\xi(+),\xi(-))$
of $Gr_0(P)$ with the loop $\xi = \xi(-) \circ \widetilde{\xi(+)}$ based
at the $*$ vertex.
Equivalently, the loop $\xi$ based at $*$ is identified with
$(\widetilde{\xi_{[n,2n]}},\xi_{[0,n]})$.
\begin{proposition}\label{intertwine}
The maps $\{\theta_n : P_{2n}(\Gamma,*) \rightarrow P_n\}_{n \geq 0}$ defined by
$$
\theta_n([\eta]) = \frac{\mu(v^\eta_0)}{\mu(v^\eta_n)} \eta
$$
for $[\eta] \in P_{2n}(\Gamma,*)$,
are $\CE_{ev}$-equivariant.
\end{proposition}

\begin{proof}
It clearly suffices to verify the intertwining assertion
on generators $S^{2n}_i$. Hence we need to check that for
$[\xi] \in P_{2n}(\Gamma,*)$ and $1 \leq i < 2n$ the equality
$$
\theta_{n-1}(S^{2n}_i([\eta])) = Z_{S^{2n}_i} (\theta_n([\eta]))
$$
holds.
There are three cases according as $i < n$, $i = n$ or $i > n$.
We will do the first case. The third is similar and the second is easier.

When $i < n$ the annular tangle $S^{2n}_i$ is shown in Figure \ref{snitangle}. The dotted lines are meant to indicate a decomposition
\begin{figure}[htbp]
\psfrag{n-1}{\tiny $n-1$}
\psfrag{i-1}{\tiny $i-1$}
\psfrag{n-i-1}{\tiny $n-i-1$}
\begin{center}
\includegraphics[height=4cm]{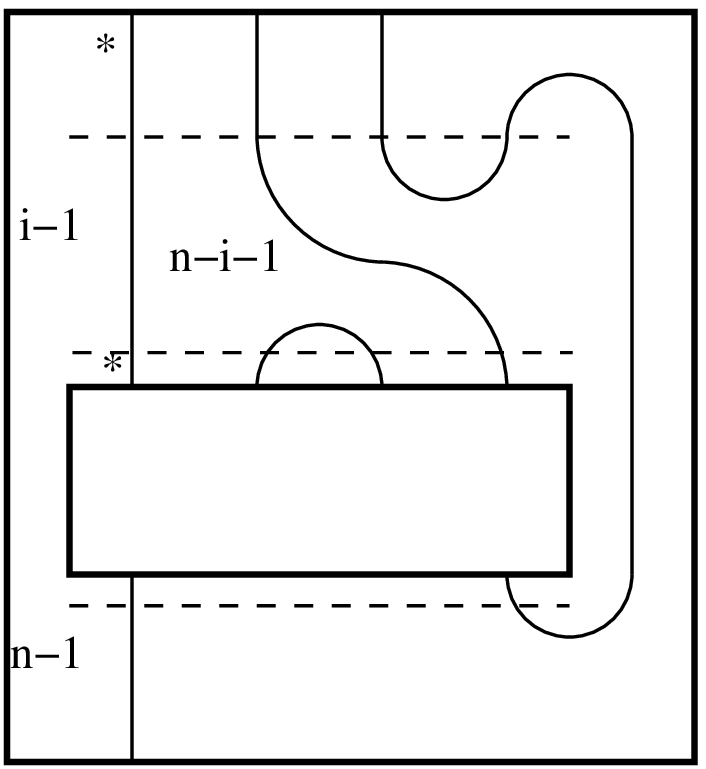}
\caption{The annular tangle $S^{2n}_i$}
\label{snitangle}
\end{center}
\end{figure}
as the (right) conditional expectation tangle applied to the
result of post-multiplication  with a Temperley-Lieb tangle. The Temperley-Lieb tangle
in question here is seen to be the product $E_{i+1}E_{i+1} \cdots E_n$,
where $E_t = \delta e_t$. It now follows by induction on $n-i$ using equations (\ref{jones}) and (\ref{inclusion}) and the multiplication in $P_n$ that
$E_{i+1}E_{i+1} \cdots E_n$ is given by
\begin{eqnarray*}
\sum_{v \in V} \sum_{\xi(\pm) \in {\mathcal P}_n(*,v)}
 \delta^{\xi(-)_{[0,i-1]}}_{\xi(+)_{[0,i-1]}}
 \delta^{\xi(-)_{[i-1,n-2]}}_{\xi(+)_{[i+1,n]}}
 \delta^{\xi(-)_{n-1}}_{\widetilde{\xi(-)_n}} 
 \delta^{\xi(+)_{i}}_{\widetilde{\xi(+)_{i+1}}} 
 \frac{\mu(v^{\xi(-)}_{n-1})\mu(v^{\xi(+)}_{i})}{\mu(v^{\xi(-)}_n) \mu(v^{\xi(+)}_{i+1})}
(\xi(+),\xi(-)).
\end{eqnarray*}
It follows that $Z_{S^{2n}_i}(\eta)$ is given by $\delta$ times the conditional expectation
onto $P_{n-1}$ of the product
\begin{eqnarray*}
\lefteqn{(\eta(+),\eta(-)) \times} \\
&\sum_{v \in V} \sum_{\xi(\pm) \in {\mathcal P}_n(*,v)}
 \delta^{\xi(-)_{[0,i-1]}}_{\xi(+)_{[0,i-1]}}
 \delta^{\xi(-)_{[i-1,n-2]}}_{\xi(+)_{[i+1,n]}}
 \delta^{\xi(-)_{n-1}}_{\widetilde{\xi(-)_n}} 
 \delta^{\xi(+)_{i}}_{\widetilde{\xi(+)_{i+1}}} 
 \frac{\mu(v^{\xi(-)}_{n-1})\mu(v^{\xi(+)}_{i})}{\mu(v^{\xi(-)}_n) \mu(v^{\xi(+)}_{i+1})}
(\xi(+),\xi(-))\\
&=
\sum_{\lambda \in {\mathcal P}_1(v^{\eta(-)}_{n},\cdot)} \delta^{\eta(-)_{i}}_{\widetilde{\eta(-)_{i+1}}}
\frac{\mu(f(\lambda))\mu(v^{\eta(-)}_{i})}{\mu(v^{\eta(-)}_n) \mu(v^{\eta(-)}_{i+1})} (\eta(+), \eta(-)_{[0,i-1]} \circ \eta(-)_{[i+1,n]} \circ \lambda \circ \widetilde{\lambda}).
\end{eqnarray*}
Now use equation (\ref{condexp}) to conclude that $Z_{S^{2n}_i}(\eta)$ is
\begin{eqnarray*}
\lefteqn{\sum_{\lambda \in {\mathcal P}_1(v^{\eta(-)}_{n},\cdot)} \delta^{\eta(-)_{i}}_{\widetilde{\eta(-)_{i+1}}}
\frac{\mu(f(\lambda))\mu(v^{\eta(-)}_{i})}{\mu(v^{\eta(-)}_n) \mu(v^{\eta(-)}_{i+1})} \delta^{\widetilde{\lambda}}_{\eta(+)_n}
\frac{\mu^2(v^{\eta(-)}_n)}{\mu^2(f(\lambda))} \times} \\
	&( \eta(+)_{[0,n-1]}, \eta(-)_{[0,i-1]} \circ \eta(-)_{[i+1,n]} \circ \lambda)\\
&=\delta^{\eta(-)_{i}}_{\widetilde{\eta(-)_{i+1}}} 
\frac{\mu(v^{\eta(-)}_n) \mu(v^{\eta(-)}_{i})}{\mu(v^{\eta(+)}_{n-1}) \mu(v^{\eta(-)}_{i+1})}
(\eta(+)_{[0,n-1]}, \eta(-)_{[0,i-1]} \circ \eta(-)_{[i+1,n]} \circ \widetilde{\eta(+)_n})\\
&=\delta^{\eta_{i}}_{\widetilde{\eta_{i+1}}} 
\frac{\mu(v^{\eta}_n) \mu(v^{\eta}_{i})}{\mu(v^{\eta}_{n+1})\mu(v^{\eta}_{i+1})}
 ~\eta_{[0,i-1]} \circ \eta_{[i+1,2n]}.
\end{eqnarray*}

Hence,
$$
Z_{S^{2n}_i}(\theta_n([\eta])) = \frac{\mu(v^\eta_0)}{\mu(v^\eta_n)} Z_{S^{2n}_i}(\eta) = \delta^{\eta_{i}}_{\widetilde{\eta_{i+1}}} 
\frac{\mu(v^{\eta}_0) \mu(v^{\eta}_{i})}{\mu(v^{\eta}_{n+1})\mu(v^{\eta}_{i+1})}
 ~\eta_{[0,i-1]} \circ \eta_{[i+1,2n]}.
$$
On the other hand, we have by definition,
$$
S^{2n}_i([\eta]) = \delta^{\eta_i}_{\widetilde{\eta_{i+1}}}
\frac{\mu(v^{\eta}_i)}{\mu(v^\eta_{i+1})} [\eta_{[0,i-1]} \circ
\eta_{[i+1,2n]}],
$$
and thus
$$
\theta_{n-1}(S^{2n}_i([\eta])) = \delta^{\eta_{i}}_{\widetilde{\eta_{i+1}}} 
\frac{\mu(v^{\eta}_0) \mu(v^{\eta}_{i})}{\mu(v^{\eta}_{n+1})\mu(v^{\eta}_{i+1})}
 ~\eta_{[0,i-1]} \circ \eta_{[i+1,2n]} = Z_{S^{2n}_i}(\theta_n([\eta])),
$$
as was to be seen.
\end{proof}

\begin{figure}[htbp]
\psfrag{1}{\tiny $1$}
\psfrag{2}{\tiny $2$}
\psfrag{3}{\tiny $3$}
\psfrag{4}{\tiny $4$}
\psfrag{5}{\tiny $5$}
\psfrag{6}{\tiny $6$}
\psfrag{7}{\tiny $7$}
\psfrag{8}{\tiny $8$}
\begin{center}
\includegraphics[height=2cm]{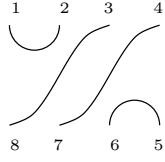}
\caption{The Temperley-Lieb tangle $T = \{ \{1,2\},\{3,8\},\{4,7\},\{5,6\}\}$}
\label{tltgl}
\end{center}
\end{figure}

Next, we generalise the path-basis expression for the Jones projections to arbitrary Temperley-Lieb tangles.
Let $T$ be a Temperley-Lieb equivalence relation on $\{1,2,\cdots,2n\}$ also identified with a Temperley-Lieb tangle as in the following example. Say $T = \{ \{1,2\},\{3,8\},\{4,7\},\{5,6\}\}$. The corresponding tangle is shown in Figure
\ref{tltgl}.
Given such a Temperley-Lieb equivalence relation $T$ we let
$T_t$ be the subset of `through classes', $T_u$ be the
subset of `up classes' and $T_d$ be the set of `down classes' of $T$, so that $T = T_t \coprod T_u \coprod T_d$. In this example, $T_t = \{\{3,8\},\{4,7\}\}$, $T_u = \{ \{1,2\}\}$ and $T_d =  \{\{5,6\}\}$.

\begin{proposition}\label{tlexp}
For any Temperley-Lieb equivalence relation $T$ on $\{1,2,\cdots,2n\}$, the element $Z_T(1) \in P_n$ is given by
\begin{eqnarray*}
\sum_{v \in V} \sum_{\xi(\pm) \in {\mathcal P}_n(*,v)}
\left( \prod_{\{i,j\} \in T_t : i<j} \delta^{\xi(-)_i}_{\xi(+)_{2n+1-j}} \right)
\left( \prod_{\{i,j\} \in T_u : i<j} \delta^{\xi(-)_i}_{\widetilde{\xi(-)_j}} \frac{\mu(v^{\xi(-)}_i)}{\mu(v^{\xi(-)}_j)} \right) \times\\
\left( \prod_{\{i,j\} \in T_d : i>j} \delta^{\xi(+)_{2n+1-i}}_{\widetilde{\xi(+)_{2n+1-j}}} \frac{\mu(v^{\xi(+)}_{2n+1-i})}{\mu(v^{\xi(+)}_{2n+1-j})} \right)
(\xi(+),\xi(-)).
\end{eqnarray*}
\end{proposition}

For instance, for the Temperley-Lieb relation $T$ of Figure \ref{tltgl},
\begin{eqnarray*}
Z_T(1) = \sum_{\xi(\pm) \in {\mathcal P}_4(*,\cdot)}
\left( \delta^{\xi(-)_3}_{\xi(+)_{1}} \delta^{\xi(-)_4}_{\xi(+)_{2}} \right)
\left( \delta^{\xi(-)_1}_{\widetilde{\xi(-)_2}} \frac{\mu(v^{\xi(-)}_1)}{\mu(v^{\xi(-)}_2)} \right)
\left( \delta^{\xi(+)_{3}}_{\widetilde{\xi(+)_{4}}} \frac{\mu(v^{\xi(+)}_{3})}{\mu(v^{\xi(+)}_{4})} \right)
(\xi(+),\xi(-)).
\end{eqnarray*}

\begin{proof}[Proof of Proposition \ref{tlexp}]
Suppose that $Z_T(1) = \sum_{v \in V} \sum_{\xi(\pm) \in {\mathcal P}_n(*,v)}
c_\xi ~ (\xi(+),\xi(-))$.
Since the $(\xi(+),\xi(-))$ are orthogonal (for the inner product on $P_n$ given by $\langle x,y \rangle = tr(y^*x)$) with $||(\xi(+),\xi(-))||^2
= \delta^{-n} \frac{\mu^2(v^{\xi(\pm)}_n)}{\mu^2(*)}$,
\begin{eqnarray*}
\overline{c_\eta} \delta^{-n} \frac{\mu^2(v^{\eta(\pm)}_n)}{\mu^2(*)} &=& \langle (\eta(+),\eta(-)), Z_T(1) \rangle\\
&=& tr(Z_T(1)^*(\eta(+),\eta(-)))\\
&=& \delta^{-n} \times {\text {picture trace of~}} Z_T(1)^*(\eta(+),\eta(-))\\
&=& \delta^{-n} \times {\text {picture trace of~}} Z_{T^*}(1)(\eta(+),\eta(-)).
\end{eqnarray*}
Hence,
\begin{equation}\label{ceta}
\overline{c_\eta} = \frac{\mu^2(*)}{\mu^2(v^{\xi(\pm)}_n)} \times {\text {picture trace of~}} Z_{T^*}(1)(\eta(+),\eta(-)).
\end{equation}

We next compute the picture trace of $Z_{T^*}(1)(\eta(+),\eta(-))$.
The equivalence relation $T^*$ is the one obtained from $T$
be replacing each $i$ by $2n+1-i$.
Regarding $T$ as an element of $Hom([2n],[0])$, there is
an associated $0_+$-annular tangle, which also we will denote by 
$T$. The context and the nature of its arguments should make it clear
whether we are referring to the morphism $T$ or the
associated Temperley-Lieb tangle $T$ or the associated
annular tangle $T$.
Some doodling now shows that 
\begin{equation}\label{pictr}
{\text {picture trace of~}} Z_{T^*}(1)(\eta(+),\eta(-)) = Z_T((\eta(+),\eta(-))) = Z_T(\eta),
\end{equation}
where, clearly, $T^*$
is regarded as a Temperley-Lieb tangle and $T$ as an annular tangle.

Finally, notice that  Proposition \ref{intertwine}
says that 
\beast \label{teta}
T \in Hom([2n],[2m]) \Rar \theta_m \circ T = T \circ \theta_n.
\eeast 
When $m=0$, since $\theta_0 = id_\C$, we see that $T = T \circ \theta_n$. i.e., 
\bea
T([\eta]) &=& T(\frac{\mu(v^\eta_0)}{\mu(v^\eta_n)} \eta)\nonumber \\
&=& \frac{\mu(v^\eta_0)}{\mu(v^\eta_n)} Z_T(\eta).
\eea

However, 
by Proposition \ref{diffexp}, we have
\begin{eqnarray*}
T([\eta]) =
\frac{\mu(v^\eta_n)}{\mu(v^\eta_0)}
  \left( \prod_{\{i,j\} \in T: i < j \leq n} \delta^{\eta_i}_{\widetilde{\eta_j}}
\frac{\mu(v^{\eta}_i)}{\mu(v^\eta_{j})} \right)
\left( \prod_{\{i,j\} \in T: i \leq n < j} \delta^{\eta_i}_{\widetilde{\eta_j}}
 \right) \times\\
   \left( \prod_{\{i,j\} \in T: n < i < j} \delta^{\eta_i}_{\widetilde{\eta_j}}
\frac{\mu(v^{\eta}_i)}{\mu(v^\eta_{j})} \right).&
\end{eqnarray*}
Putting this together with equations (\ref{ceta}),(\ref{pictr}) and (\ref{teta}) yields
$$
\overline{c_\eta} = \left( \prod_{\{i,j\} \in T: i < j \leq n} \delta^{\eta_i}_{\widetilde{\eta_j}}
\frac{\mu(v^{\eta}_i)}{\mu(v^\eta_{j})} \right)
\left( \prod_{\{i,j\} \in T: i \leq n < j} \delta^{\eta_i}_{\widetilde{\eta_j}}
 \right)
   \left( \prod_{\{i,j\} \in T: n < i < j} \delta^{\eta_i}_{\widetilde{\eta_j}}
\frac{\mu(v^{\eta}_i)}{\mu(v^\eta_{j})} \right).
$$
Observing that $c_\eta$ is real and comparing with the statement of the proposition finishes the proof.
(Note that when $\{i,j\} \in T$ with $n < i < j$,
$$
\delta^{\eta_i}_{\widetilde{\eta_j}}
\frac{\mu(v^{\eta}_i)}{\mu(v^\eta_{j})} =
\delta^{\eta_i}_{\widetilde{\eta_j}}
\frac{\mu(v^{\eta}_{j-1})}{\mu(v^\eta_{i-1})} =
\delta_{\xi(+)_{2n+1-j}}^{\widetilde{\xi(+)_{2n+1-i}}} \frac{\mu(v^{\xi(+)}_{2n+1-j})}{\mu(v^{\xi(+)}_{2n+1-i})},
$$
which is to be compared with the third product term in the statement.)
\end{proof}

We will now write the structure maps of the algebra $Gr_0(P)$
of \cite{GnnJnsShl2008} in terms of the path bases for the $P_n$.
Recall that the algebra $Gr_0(P) = \oplus_{n=0}^\infty P_n$ is
a graded algebra with the multiplication map $\bullet : P_m \otimes P_n \rightarrow P_{m+n}$ given by the tangle in Figure \ref{mult} below.
\begin{figure}[htbp]
\psfrag{2m}{\tiny $2m$}
\psfrag{2n}{\tiny $2n$}
\begin{center}
\includegraphics[height=2cm]{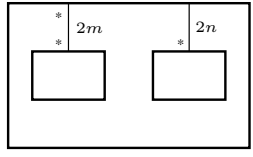}
\caption{Multiplication in $Gr_0(P)$}
\label{mult}
\end{center}
\end{figure}

\begin{proposition}\label{multiplicativity}
For paths $\xi \in P_m \subseteq Gr_0(P)$ and $\eta \in P_n \subseteq Gr_0(P)$,
$$
\xi \bullet \eta = \frac{\mu(v^\xi_m)\mu(v^\eta_n)}{\mu(v^{\xi \circ \eta}_{m+n}) \mu(v^\eta_{0})}\  \xi \circ \eta 
$$
\end{proposition}

\begin{proof}
We will deal with the case $m \geq n$. The other case is similar.
The tangle of Figure \ref{mult} can be expressed in
terms of the inclusion, Temperley-Lieb and multiplication tangles
as in Figure~\ref{multexp}.
\begin{figure}[htbp]
\psfrag{m}{\tiny $m$}
\psfrag{n}{\tiny $n$}
\psfrag{m-n}{\tiny $m-n$}
\begin{center}
\includegraphics[height=4cm]{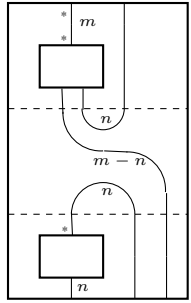}
\caption{Standard tangle expression of tangle in Figure \ref{mult}}
\label{multexp}
\end{center}
\end{figure}
Recall that in a tangle picture, a non-negative integer $t$ written
beside a string indicates a $t$-cable of the string.
We see from this figure that the product of $\xi$ and $\eta$ in
$Gr_0(P)$ is a product of three terms in $P_{m+n}$, namely, $\xi$
included into $P_{m+n}$, a Temperley-Lieb tangle and $\eta$
included into $P_{m+n}$.

It now follows from Proposition \ref{tlexp} and equation (\ref{inclusion}) that
\begin{eqnarray*}
\xi \bullet \eta &=&  \left( \sum_{\rho \in {\mathcal P}_m(f(\eta(\pm)),\cdot)} (\eta(+) \circ \rho,\eta(-) \circ \rho)\right) \times\\
& & \left\{\sum_{v \in V} \sum_{\zeta(\pm) \in {\mathcal P}_{m+n}(*,v)}
\left( \delta^{\zeta(-)_{[0,m-n]}}_{\zeta(+)_{[2n,m+n]}} \right)\right. \times \\
& &\left. \left( \delta^{\zeta(-)^{[m-n,m]}}_{\widetilde{\zeta(-)_{[m,m+n]}}} \frac{\mu(v^{\zeta(-)}_m)}{\mu(v^{\zeta(-)}_{m+n})} \right)
\left( \delta^{\zeta(+)_{[0,n]}}_{\widetilde{\zeta(+)_{[n,2n]}}} \frac{\mu(v^{\zeta(+)}_{n})}{\mu(v^{\zeta(+)}_{2n})} \right)
(\zeta(+),\zeta(-))\right\} \times\\
& & \left( \sum_{\lambda \in {\mathcal P}_n(f(\xi(\pm)),\cdot)} (\xi(+) \circ \lambda,\xi(-) \circ \lambda)\right).
\end{eqnarray*}

Since the path basis elements multiply as matrix units, given  $\lambda,\zeta(\pm)$ and $\rho$, the product
of the terms corresponding to these in the above expression is
non-zero only if the following equations hold.
\begin{eqnarray*}
\zeta(-) &=& \xi(+) \circ \lambda\\
\zeta(+) &=& \eta(-) \circ \rho\\
\zeta(-)_{[0,m-n]} &=& \zeta(+)_{[2n,m+n]}\\
\zeta(-)_{[m-n,m]} &=& \widetilde{\zeta(-)_{[m,m+n]}}\\
\zeta(+)_{[0,n]} &=& \widetilde{\zeta(+)_{[n,2n]}}
\end{eqnarray*}

A little thought now shows that the following equations are
consequences.
\begin{eqnarray*}
\zeta(-)_{[0,m]} &=& \xi(+)\\
\zeta(+)_{[0,n]} &=& \eta(-)\\
\zeta(-)_{[m,m+n]} &=& \widetilde{\zeta(-)_{[m-n,m]}} = \widetilde{\xi(+)_{[m-n,m]}}\\
\zeta(+)_{[n,2n]} &=& \widetilde{\zeta(+)_{[0,n]}} = \widetilde{\eta(-)}\\
\zeta(+)_{[2n,m+n]} &=& \zeta(-)_{[0,m-n]} = \xi(+)_{[0,m-n]}\\
\lambda &=& \zeta(-)_{[m,m+n]} = \widetilde{\xi(+)_{[m-n,m]}}\\
\rho &=& \zeta(+)_{[n,m+n]} = \widetilde{\eta(-)} \circ \xi(+)_{[0,m-n]}
\end{eqnarray*}

Thus, exactly one term is non-zero, which corresponds to
$\lambda = \widetilde{\xi(+)_{[m-n,m]}}$, $\rho = \widetilde{\eta(-)} \circ \xi(+)_{[0,m-n]}$, $\zeta(-) = \xi(+) \circ \widetilde{\xi(+)_{[m-n,m]}}$
and $\zeta(+) = \eta(-) \circ \widetilde{\eta(-)} \circ \xi(+)_{[0,m-n]}$.
Hence
$$
\xi \bullet \eta = \frac{\mu(v^{\xi(+)}_m)}{\mu(v^{\xi(+)}_{m+n})} \frac{\mu(v^{\eta(-)}_n)}{\mu(v^{\eta(-)}_0)} (\eta(+) \circ \widetilde{\eta(-)} \circ \xi(+)_{[0,m-n]}, \xi(-) \circ \widetilde{\xi(+)_{[m-n,m]}})
$$
Noting now that $\xi = \xi(-) \circ \widetilde{\xi(+)}$, $\eta = \eta(-) \circ \widetilde{\eta(+)}$ and $\xi \circ \eta = \xi(-) \circ \widetilde{\xi(+)} \circ \eta(-) \circ \widetilde{\eta(+)}$, the proof is finished.
\end{proof}

\begin{proposition}\label{griso}
The map $\theta : Gr(\Gamma,*) \rightarrow Gr_0(P)$ defined for
$[\xi] \in P_{2n}(\Gamma,*)$ by $$\theta([\xi]) = \frac{\mu(v^\xi_0)}{\mu(v^\xi_n)} \xi \in P_n$$ and extended linearly
is an isomorphism of graded, $*$-probability spaces.
\end{proposition}

\begin{proof}
That $\theta$ is a graded, linear isomorphism is clear.
Multiplicativity of $\theta$ follows from Proposition \ref{multiplicativity}
while $*$-preservation is straightforward. To verify that 
$\theta$ preserves trace, note that by definition of the trace $\tau$ in
$Gr(\Gamma,*)$, for $[\xi] \in P_{2n}(\Gamma,*)$,
$$
\tau([\xi]) =  \sum_T \tau_T([\xi]),
$$
where the sum is over all Temperley-Lieb equivalence
relations
$T$ on $\{1,2,\cdots,2n\}$ and $\tau_T([\xi])$ is (from the proof of
Proposition \ref{same}) $\frac{1}{\mu^2(*)} t \circ T([\xi])$ where $t : P_0(\Gamma) \rightarrow {\mathbb C}$ is the linear extension of the
map taking $[(v)]$ to $\mu^2(v)$. Identifying $P_0(\Gamma,*)$ with
${\mathbb C}$, $\tau_T([\xi]) = T([\xi]).$
Equations (\ref{pictr}) and (\ref{teta}) now imply that $\tau_T([\xi])$ is 
$\frac{\mu(v^\xi_0)}{\mu(v^\xi_n)}$ times the picture trace of $Z_{T^*}(1)\xi$. Summing over all Temperley-Lieb equivalence
relations gives by definition the trace of
$\frac{\mu(v^\xi_0)}{\mu(v^\xi_n)}\xi$ in $Gr_0(P)$, as desired.
\end{proof}

We  apply this proposition and Theorem \ref{mgamlf} to the GJS construction.

\begin{theorem}\label{main}
Let $P$ be a subfactor planar algebra of finite depth and modulus $\delta > 1$, and $M_0$ be the factor
constructed from $P$ by the construction in \cite{GnnJnsShl2008}.
Then, $M_0 \cong LF(r)$ for some $1 < r < \infty$.
\end{theorem}

\begin{proof}
Let $\Gamma$ be the (finite) principal graph of $P$ equipped with the
Perron-Frobenius weighting, so that by Theorem \ref{mgamlf}, $M(\Gamma)$
is $LF(t)$ for some $1 < t < \infty$. Now Proposition \ref{griso}
implies that $M(\Gamma,*)$ is isomorphic to $M_0$ and
so
$M_0 \cong LF(r)$ for some $1 < r < \infty$.
\end{proof} 

Our final result is an analogue of Theorem \ref{main} for the factor
$M_1$ constructed from $P$.

\begin{theorem}\label{main2}
Let $P$ be a subfactor planar algebra of finite depth and modulus $\delta > 1$, and $M_0 \subseteq M_1$ be the subfactor
constructed from $P$ by the construction in \cite{GnnJnsShl2008}.
Then $M_1 \cong LF(s)$ for some $1 < s < \infty$.
\end{theorem}

Since the proof is very similar to that of Theorem \ref{main}, we will only sketch the proof  giving details where it differs from the previous
proof.  We first recall some preliminaries from \cite{KdySnd2004}.

There is an `operation on tangles' denoted $T \mapsto T^{-}$ which
moves the $*$-region of each of its boxes anti-clockwise by 1 and reverses the shading. There is an associated `operation on planar algebras' denoted $P \mapsto {^-P}$ where $^-P$ is the planar
algebra with spaces
\begin{eqnarray*}
^{-}P_{0_\pm} &=& P_{0_\mp}\\
^{-}P_k &=& P_k, ~~~~~k > 0,
\end{eqnarray*}
and tangle action defined by $Z^{^-P}_T = Z^P_{T^-}$.
If $P$ is a subfactor planar algebra, then so is $Q = {^-P}$ and further
$^-Q$ is isomorphic to $P$.

\begin{figure}[htbp]
\psfrag{T}{\tiny $T$}
\psfrag{2n}{\tiny $2n-2$}
\begin{center}
\includegraphics[height=2.5cm]{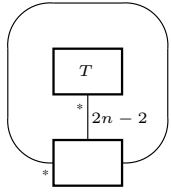}
\caption{$T$-component of the trace  in $^-Gr_1(P)$}
\label{gr1tr}
\end{center}
\end{figure}

Now, given a subfactor planar algebra $P$, we define a graded, non-commutative probability space $^-Gr_1(P)$ associated to
$P$ as follows. As a vector space $^-Gr_1(P) = \oplus_{n \geq 1} P_n$. The multiplication map $\bullet : P_m \otimes P_n \rightarrow P_{m+n-1}$ is defined by the tangle in Figure \ref{gr1mult} below.
\begin{figure}[htbp]
\psfrag{2m}{\tiny $2m-2$}
\psfrag{2n}{\tiny $2n-2$}
\begin{center}
\includegraphics[height=2cm]{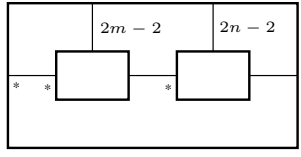}
\caption{Multiplication in $^-Gr_1(P)$}
\label{gr1mult}
\end{center}
\end{figure}
The adjunction map in $^-Gr_1(P)$ restricts to the adjunction maps
in $P_n$ for each $n \geq 1$. A trace is defined in $^-Gr_1(P)$ by
letting $\tau(\xi)$ for $\xi \in P_n \subseteq {^-Gr_1(P)}$ be the
sum over all Temperley-Lieb tangles $T$ of the scalar defined by Figure \ref{gr1tr}.

The structure maps of $^-Gr_1(P)$ are all derived from those of
$Gr_1(P)$ (see \cite{GnnJnsShl2008}) using the operation $^-$.
Observe that the vector space underlying both $^-Gr_1(P)$ and $Gr_1({^-P})$ is the same, namely, $\oplus_{n \geq 1} P_n$.
A little thought now yields the following.

\begin{lemma}\label{bargr}
For any subfactor planar algebra $P$, the tracial $*$-probability
spaces $^-Gr_1(P)$ and $Gr_1({^-P})$ are isomorphic
by the identity map of the underlying vector spaces.\qed
\end{lemma}

Applying Lemma \ref{bargr} with $Q = {^-P}$ in place of $P$ and
using that ${^-Q} \cong P$ shows that $^-Gr_1(Q) \cong Gr_1({P})$
as probability spaces. We now proceed towards an analogue of
Proposition \ref{griso} for $^-Gr_1(Q)$. Let $\overline{\Gamma}$
denote the principal graph of $Q$. Since $P$ is of finite depth, so is $Q$, and thus $\overline{\Gamma}$ is a finite graph. Equip $\overline{\Gamma}$ with its Perron-Frobenius
weighting.  
A basis of $Q_n$ is then given by pairs of paths $(\xi(+),\xi(-))$ in $\overline{\Gamma}$
such that 
$\xi(\pm)$ are paths of length $n$ in $\overline{\Gamma}$
beginning at its $*$ and having the same end-point.
Again, we identify the basis element $(\xi(+),\xi(-))$ with the loop
$\xi(-) \circ \widetilde{\xi(+)}$ based at $*$.

Observe that the $0^{th}$-graded piece of $^-Gr_1(Q)$ can be identified as an algebra with $Q_1^{op}$. In particular, each vertex
$v$ in $\overline{\Gamma}$ at distance 1 from its $*$ vertex
gives a minimal central projection $f(v,1)$ in $[^-Gr_1(Q)]_1 = Q_1^{op}$
and we denote by $q(v,1)$ any minimal projection of $Q_1^{op}$
lying under $f(v,1)$. A choice of
$q(v,1)$ is the matrix unit $(\nu, \nu) \in Q_1^{op}$ where
$\nu$ is any path of length 1 in $\overline{\Gamma}$ from
$*$ to $v$. We fix this choice.

The following proposition expresses the multiplication of  $^-Gr_1(Q)$ in terms of its path basis.

\begin{proposition}\label{gr1multiplicativity}
For paths $\xi \in Q_m \subseteq {^-Gr}_1(Q)$ and $\eta \in P_n \subseteq {^-Gr}_1(Q)$,
$$
\xi \bullet \eta = \delta^{\xi_{2m}}_{\widetilde{\eta_1}}\frac{\mu(v^\xi_m)\mu(v^\eta_n)}{\mu(v^{\xi \circ \eta}_{m+n-1}) \mu(v^\eta_{1})}\  \xi_{[0,2m-1]} \circ \eta _{[1,2n]}
$$
\end{proposition}

\begin{proof}[Sketch of proof] Suppose that $m \geq n$. The key fact
is that the tangle of Figure \ref{gr1mult} is expressible in terms
of the inclusion, Temperley-Lieb and multiplication tangles as in
Figure \ref{gr1multexp}.
\begin{figure}[htbp]
\psfrag{m}{\tiny $m$}
\psfrag{n}{\tiny $n-1$}
\psfrag{n+1}{\tiny $n$}
\psfrag{m-n}{\tiny $m-n$}
\begin{center}
\includegraphics[height=4cm]{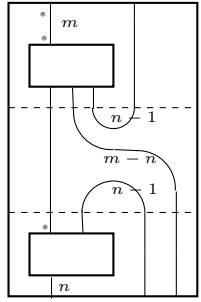}
\caption{Multiplication of $^-Gr_1(Q)$ in terms of standard tangles}
\label{gr1multexp}
\end{center}
\end{figure}
We omit the rest of the proof which is very similar to that of Proposition \ref{multiplicativity}.
\end{proof}

It follows from Proposition \ref{gr1multiplicativity} that a basis of
$q(v,1)(^-Gr_1(Q))q(v,1)$ is given by all paths of the form $\nu \circ
\xi \circ \widetilde{\nu}$ where $\xi$ ranges over all paths in $\overline{\Gamma}$ from $v$ to $v$. This is suggestive of the following key
isomorphism which is the analogue of Proposition~\ref{griso} and whose proof is similar (and omitted).

\begin{proposition}\label{gr1iso}
The map $\theta : Gr(\overline{\Gamma},v) \rightarrow q(v,1)(^-Gr_1(Q))q(v,1)$ defined for
$[\xi] \in P_{2n}(\overline{\Gamma},v)$ by 
$$
\theta([\xi]) = \frac{\mu(v^\xi_0)}{\mu(v^\xi_n)} \nu \circ \xi \circ \widetilde{\nu} \in Q_{n+1}
$$ 
and extended linearly
is an isomorphism of graded, $*$-probability spaces. \qed
\end{proposition}

We conclude with the proof of Theorem \ref{main2}.

\begin{proof}[Proof of Theorem \ref{main2}]
From Proposition \ref{gr1iso} and the isomorphism of ${^-Gr}_1(Q)$
with $Gr_1(P)$, it follows by completing that some corner of
$M_1$ is isomorphic to $M(\overline{\Gamma},v)$ - a corner
of $M(\overline{\Gamma})$.
Since $M(\overline{\Gamma})$ is an interpolated free group factor
with finite parameter 
by Theorem \ref{mgamlf}, the proof is complete.
\end{proof}

\section*{Acknowledgment}
We are deeply indebted to Ken Dykema for his constant guidance
throughout the preparation of this manuscript, and even more
for working overtime to identify and prove the statement needed
in \S 5.

\end{document}